\input amstex
\magnification 1200
\TagsOnRight
\def\qed{\ifhmode\unskip\nobreak\fi\ifmmode\ifinner\else
 \hskip5pt\fi\fi\hbox{\hskip5pt\vrule width4pt
 height6pt depth1.5pt\hskip1pt}}
\NoBlackBoxes
\baselineskip 17.9 pt
\parskip 5 pt
\def\stretch {\noalign{\medskip}}
\define \bC {\bold C}
\define \bCp {\bold C^+}
\define \bCm {\bold C^-}
\define \bCpb {\overline{\bold C^+}}
\define \ds {\displaystyle}
\define \bR {\bold R}
\define \bm {\bmatrix}
\define \endbm {\endbmatrix}

\centerline {\bf INVERSE PROBLEM WITH TRANSMISSION EIGENVALUES}
\vskip -5 pt
\centerline {\bf FOR THE DISCRETE SCHR\"ODINGER EQUATION}

\vskip 7 pt
\centerline {Tuncay Aktosun}
\vskip -8 pt
\centerline {Department of Mathematics}
\vskip -8 pt
\centerline {University of Texas at Arlington}
\vskip -8 pt
\centerline {Arlington, TX 76019-0408, USA}
\vskip -8 pt
\centerline {aktosun\@uta.edu}

\centerline {Vassilis G. Papanicolaou}
\vskip -8 pt
\centerline {Department of Mathematics}
\vskip -8 pt
\centerline {National Technical University of Athens}
\vskip -8 pt
\centerline {Zografou Campus}
\vskip -8 pt
\centerline {157 80, Athens, Greece}
\vskip -8 pt
\centerline {papanico\@math.ntua.gr}

\noindent {\bf Abstract}:
The discrete Schr\"odinger equation
with the Dirichlet boundary condition is considered on a
half-line lattice when the potential is real valued and compactly supported.
The inverse problem of recovery of
the potential from the so-called transmission eigenvalues is
analyzed. The Marchenko method and the Gel'fand-Levitan method
are used to solve the inverse problem uniquely, except in one
``unusual" case where the sum of the transmission eigenvalues
is equal to a certain integer related to the support of the potential.
It is shown that in the unusual case there may be a unique solution
corresponding to certain sets of transmission eigenvalues, there may be a
finite number of distinct potentials for some sets of transmission eigenvalues, or there may be infinitely many potentials
for some sets of transmission eigenvalues. The theory presented
is illustrated with several explicit examples.

\vskip -8 pt
\par\noindent {\bf Short title:} Inverse problem with transmission eigenvalues
\vskip -8 pt
\par\noindent {\bf Keywords:} Discrete Schr\"odinger equation,
transmission eigenvalues, inverse problem,
Marchenko method, Gel'fand-Levitan method, spectral function

\newpage

\noindent {\bf 1. INTRODUCTION}
\vskip 3 pt

Let us consider the discrete Schr\"odinger equation on the half line
$$-\psi_{n+1}+2 \psi_n-\psi_{n-1}+V_n \psi_n=\lambda\, \psi_n,\qquad n\ge 1,\tag 1.1$$
where $\lambda$ is the spectral parameter, $n$ is the spacial independent variable
taking positive integer values, and
the subscripts are used to denote the dependence on $n.$
Thus, $\psi_n$ denotes the wavefunction
$\psi(\lambda,n).$ Note that $\psi_0$ appears in (1.1) when $n=1,$ and hence the value at $n=0$ of
the wavefunction
$\psi_n$ must be compatible with (1.1).
By $V_n,$ which can also be written as $V(n),$ we denote the value
of the potential at $n.$ The potential $V$ is assumed to belong to
class $\Cal A_b,$ which is specified in the following definition.

\noindent {\bf Definition 1.1} {\it The potential $V$
belongs to class $\Cal A_b$ if the $V_n$-values are real and the support
of the
potential $V$ is confined to the discrete set $n\in\{1,2,\dots,b\}$ for some positive integer $b,$ i.e. $V_n=0$ for $n>b.$}

As indicated in Definition~1.1, we assume that the potential in
(1.1) is compactly supported, i.e. it vanishes for $n\ge b+1.$
Thus, the knowledge of
the potential $V$ in class $\Cal A_b$ is equivalent to the knowledge of the ordered set
$\{V_1,V_2,\dots,V_b\}.$ We do not necessarily assume that all the values in this ordered set
are nonzero.

Note that (1.1) is the discrete analog of
the Schr\"odinger equation on the half line
$$-\psi''(k,x)+V(x)\,\psi(k,x)=k^2\,\psi(k,x), \qquad x>0,
\tag 1.2$$
where $\lambda:=k^2,$
the prime denotes the $x$-derivative, the potential
$V$ is real valued, and the support of the potential is confined to the
interval $x\in[0,b].$

We are interested in the inverse problem of recovery of
the potential $V$ in (1.1)
when the available input
data set consists of the so-called transmission eigenvalues. The precise meaning
of transmission eigenvalues is given in Section~5.
The transmission eigenvalues are uniquely determined by the potential
$V$ and the Dirichlet boundary condition (2.2) specified in
Section~2.
Informally speaking,
the transmission eigenvalues correspond to certain $\lambda$-values at which the scattering
 from (1.1) with the Dirichlet
 boundary condition (2.2)
 agrees with the scattering when the potential is identically zero.
This inverse problem was recently analyzed in [22], and it was shown that
the transmission eigenvalues uniquely determine the potential with the exception of
one case  [22], which we now call the ``unusual case,"
namely when the sum of the transmission eigenvalues is equal to $4(b-1),$ where $b$ is the positive integer appearing in Definition~1.1 related to
the support of the potential.

A recursive procedure was presented
in [22] to determine the potential from the transmission eigenvalues
in the ``usual case," and the nonuniqueness in the unusual case
was illustrated via an example in [22].
However, the procedure of [22] to
construct the potential is not as direct as either of the
Marchenko procedure and the Gel'fand-Levitan procedure
we present in this paper because it requires first the construction of the right-hand side
of (2.70) in our paper, namely the regular solution
to (1.1) as a function of $\lambda,$ and then using
that regular solution in (1.1) to extract the
values of the potential.

We remark that, when the potential $V$ vanishes for $n\ge b+1$ while
$V_b\ne 0,$ as indicated in Theorem~5.2,
the number of
transmission eigenvalues (including multiplicities) is $2b-2,$ and hence the inverse problem dealing with transmission
eigenvalues is meaningful only when $b\ge 2.$

One of our main goals in this paper is to develop a comprehensive approach to solve the
inverse problem in the usual case in order to determine
the potential $V$ in class $\Cal A_b$
 from the set of corresponding transmission eigenvalues. This is done
via the Marchenko method and also the Gel'fand-Levitan method for (1.1).
Another goal we have is
to elaborate on the unusual case and to show that the inverse problem
in the unusual case may or may not have a unique solution.
This is done by presenting various examples in Section~6 and
by showing that all of
the following possibilities occur. There may be only one real-valued
potential
corresponding to a given set of transmission eigenvalues,
there may be a finite number of distinct real-valued potentials
corresponding to a given set of transmission eigenvalues,
or there may be infinitely many distinct real-valued potentials
corresponding to a given set of transmission eigenvalues.

The literature on transmission eigenvalues and related inverse
problems is expanding rapidly. We refer the reader to [1,3,4,7,8,12-15,19-21] and the references
therein for
a historical account of transmission eigenvalues
and some important developments
in the field.
Since we are dealing with the difference equation (1.1) rather than
the differential equation (1.2), the most relevant reference
for our paper is [22]. The techniques used in our paper
are analogs of the techniques from [3,4] used in the
continuous case for (1.2).

For the inverse problem
for the discrete Schr\"odinger equation and the related
Gel'fand-Levitan method, we refer the reader to
[9], even though the discrete Schr\"odinger
equation used in [9]
differs from (1.1). In [9], the discrete
Schr\"odinger equation is treated as the Jacobi equation
$$\ds\frac{1}{2}\, e^{-(V_n+V_{n+1})/2}\phi_{n+1}+\ds\frac{1}{2}\, e^{-(V_{n-1}+V_n)/2}
\phi_{n-1}
=\lambda\,\phi_n, \qquad n\ge 1.\tag 1.3$$
To help the reader, in presenting the
results for (1.1) we at times state the corresponding results for
(1.2) by illustrating the similarities and differences between the
discrete and continuous cases.

Our paper is organized as follows. In Section~2 the preliminaries
are presented for the Marchenko method to be given in Section~3,
the Gel'fand-Levitan method to be given in Section~4, and
the inverse problem with transmission eigenvalues to be given
in Section~5. This is done by introducing the so-called
Jost solution $f_n$ to (1.1),
the regular solution $\varphi_n$ to (1.1), the Jost function $f_0$
corresponding to the value of the Jost solution at $n=0,$
and the scattering matrix $S$ given in (2.45),
and by presenting certain relevant properties of
such quantities because those properties are needed in later sections.
Note that, as in the continuous case, we call $S$ a ``matrix"
even though it is scalar valued.
In Section~3 the linear system of Marchenko equations for (1.1) with the boundary
condition (2.2) is derived and it is shown how the potential
can be recovered from the solution to the Marchenko system (3.12).
The Marchenko system uses as input the scattering matrix and
the bound-state information consisting of bound-state
energies and the so-called Marchenko bound-state norming constants.
In Section~3 it is also shown that, as a result of the
compact-support property of the potential,
the bound-state information is actually contained in the
scattering matrix and hence the scattering matrix alone
uniquely determines the potential.
In Section~4 the Gel'fand-Levitan method for (1.1)
with the boundary
condition (2.2) is derived and it is shown how the potential
can be recovered from the solution to the linear system of
Gel'fand-Levitan equations (4.17).
It is also shown how the kernel of the Gel'fand-Levitan system
can be constructed from
the Jost function alone.
In Section~5 the transmission eigenvalues related to
(1.1) with the boundary
condition (2.2) are introduced and they are shown to
correspond to the zeros of the key quantity
$D$ given in (5.2) or equivalently to the zeros of
the quantity $E$ defined in (5.7). It is also shown how the
Jost function $f_0$ can be recovered from the key quantity
$E$ in the usual case, and hence, having $f_0$ at hand,
one can use either the Marchenko method or the Gel'fand-Levitan method
to recover the potential. In Section~5 it is further
indicated what happens when the unusual case occurs,
which may prevent the recovery of a unique potential
 from the set of transmission eigenvalues.
Finally, in Section~6 some examples are presented to
illustrate the unique recovery of the potential from
the transmission eigenvalues via the Marchenko
method and the Gel'fand-Levitan method in the
usual case. In that section some further examples
are presented
in the unusual case to illustrate
that there may exist a unique real-valued potential
for some set of transmission eigenvalues, there may be
a finite number of distinct real-valued potentials for some
other sets of transmission eigenvalues, or there may be infinitely
many distinct real-valued potentials corresponding to some sets of
transmission eigenvalues.

\vskip 5 pt
\noindent {\bf 2. PRELIMINARIES}
\vskip 3 pt

In this section we introduce the Jost solution $f_n$
to the discrete Schr\"odinger
equation (1.1), the regular solution $\varphi_n$ to (1.1), the corresponding Jost function
$f_0,$
the related scattering matrix $S,$ and we present some properties of these
quantities needed later on for the solution to the inverse problem of recovery of
the potential in (1.1) from the so-called transmission eigenvalues.

Using the shift operators $P^+$ and $P^-$ and the identity operator $I,$ which are
defined as
$$P^+\psi_n=\psi_{n+1},\quad P^-\psi_n=\psi_{n-1},\quad I\psi_n=\psi_n,$$
let us write (1.1) in the operator form as $\Cal L\psi_n=\lambda \psi_n,$ where we have defined
$$\Cal L:=-P^++2 I-P^- +V_n.\tag 2.1$$
Let us associate
with (1.1) the Dirichlet boundary condition
$$\psi_0=0,\tag 2.2$$
which is the analog of the Dirichlet boundary condition
for (1.2), namely
$$\psi(0)=0.\tag 2.3$$

Note that $\Cal L$
given in (2.1) with the boundary condition (2.2) is a bounded operator
because the shift operators $P^+$ and $P^-$ appearing in
(2.1) each have norm one and the multiplication
by $V_n$ is a bounded operator. Furthermore, $\Cal L$ is
selfadjoint [23],
which can be verified by establishing the equality
$$\ds\sum_{n=1}^\infty \phi_n^\ast\,(\Cal L\,\psi_n)
=\ds\sum_{n=1}^\infty (\Cal L\,\phi_n)^\ast\,\psi_n,$$
where $\{\psi_n\}_{n=1}^\infty$ and $\{\phi_n\}_{n=1}^\infty$ are two square-summable
sequences satisfying (2.2). Note that we use an asterisk
to denote complex conjugation, and we recall that
the square summability of $\psi_n$
is described by
$$\sum_{n=1}^\infty |\psi_n|^2<+\infty.\tag 2.4$$
The bound states of $\Cal L$
correspond to the discrete spectrum of $\Cal L,$ i.e. the
set of those $\lambda$-values
for which there exists a corresponding square-summable solution to
(1.1) satisfying (2.2).

When the potential $V$ belongs to
class $\Cal A_b$ specified in Definition~1.1, as we will
 see in Theorem~2.5, the bound states can only occur
when $\lambda\in (-\infty,0)\cup (4,+\infty)$ and the number of
bound states is either zero or
a finite number. We will use $N$ to denote the number
of bound states and assume that the bound states occur
at $\lambda=\mu_j$ for $j=1,\dots,N.$
The continuous spectrum of the discrete
Schr\"odinger operator $\Cal L$ with the Dirichlet boundary condition
(2.2) is the interval $\lambda\in[0,4]$
because as we will see, for each $\lambda\in[0,4],$ we have a solution to
(1.1) satisfying (2.2) in such a way that such a solution
is either bounded in $n$ without being square summable or it
grows linearly in $n$ as $n\to+\infty.$

In case $V_n\equiv 0$ in (1.1), let us use $\overset{\circ}\to{\psi}_n$ to denote
the corresponding wavefunction. Thus, $\overset{\circ}\to{\psi}_n$ satisfies the
unperturbed discrete Schr\"odinger equation
$$-\overset{\circ}\to\psi_{n+1}+2 \overset{\circ}\to\psi_n-\overset{\circ}\to\psi_{n-1}=\lambda \overset{\circ}\to\psi_n,\qquad n\ge 1.\tag 2.5$$

It is convenient at times to use not $\lambda$ but another
spectral parameter related to $\lambda,$ which is usually
denoted by $z.$ This is done by
letting
$$z:=1-\ds\frac{\lambda}{2}+\ds\frac{1}{2}\,\sqrt{\lambda(\lambda-4)},\tag 2.6$$
where the square root is used to denote the principal
branch of the complex
square-root function.
Let us use $\bold T$ to denote the unit
circle $|z|=1,$ use $\bold T^+$ for the upper
portion of $\bold T$ given by $z=e^{i\theta}$ with
$\theta\in(0,\pi),$ and
use $\overline{\bold T^+}$ for the closure
of $\bold T^+$ given by $z=e^{i\theta}$ with
$\theta\in[0,\pi].$
Under the transformation $\lambda\mapsto z,$
the point $\lambda=\infty$
corresponds to
$z=0$ in the complex $z$-plane,
the point $\lambda=0$ corresponds to $z=1,$ and the point $\lambda=4$
to $z=-1.$
The real half line $\lambda\in(-\infty,0)$ is mapped to the real interval
$z\in (0,1)$ in the complex $z$-plane, the real interval $\lambda
\in(4,+\infty)$ is mapped to the real interval $z\in(-1,0),$ and the real interval
$\lambda\in(0,4)$ is mapped to $\bold T^+.$
Hence, the real $\lambda$-axis
is mapped to the
boundary of the upper half
of the unit disc in the complex $z$-plane.
Thus, the bound states in the complex $z$-plane can only occur
in the interval $z\in(-1,0)$ or $z\in(0,1).$
 From (2.6) it follows that
 $$z^{-1}=1-\ds\frac{\lambda}{2}-\ds\frac{1}{2}\,\sqrt{\lambda(\lambda-4)}.\tag 2.7$$
We see from (2.6) and (2.7) that
$$\lambda=2-z-z^{-1}.\tag 2.8$$

Let us express $z$ in terms of still another spectral parameter $\theta$ as
$$z=e^{i\theta}.\tag 2.9$$
Letting
$$\cos \theta=1-\ds\frac{\lambda}{2},\quad
\sin\theta=\ds\frac{1}{2}\,\sqrt{\lambda(4-\lambda)},\tag 2.10$$
as we have already observed, as $\theta$ takes values in the interval
$\theta\in[0,\pi]$ the spectral parameter $z$ traces $\overline{\bold T^+}$
and from (2.10) we see that the spectral
parameter $\lambda$ then moves in the interval $\lambda\in[0,4]$ from the left endpoint
$\lambda=0$ to the right endpoint $\lambda=4.$
Note that, from (2.8) and (2.9), we have
$$\lambda=2-2\cos\theta.\tag 2.11$$

Using (2.8) it is convenient to write (1.1) as
$$\psi_{n+1}+\psi_{n-1}=(z+z^{-1}+V_n)\,\psi_n,\qquad n\ge 1.\tag 2.12$$
Thus, (2.5) can be written as
$$\overset{\circ}\to\psi_{n+1} +\overset{\circ}\to\psi_{n-1}=(z+z^{-1})
\overset{\circ}\to\psi_n,\qquad n\ge 1.\tag 2.13$$
One can directly verify that
$z^n$ and $z^{-n}$ are solutions to
(2.13) and that they are linearly independent when $z\ne \pm 1.$ Hence,
the general solution to (2.13) can be written as a linear combination
of $z^n$ and $z^{-n}$ when $z\ne \pm 1.$ In (2.52) we display
two linearly independent solutions to (2.13) when $z=1.$
In (2.56) we list two linearly independent solutions to (2.13)
when $z=-1.$

There are certain relevant solutions to (2.12) and
hence equivalently to (1.1). One such solution is
the so-called regular solution $\varphi_n$ satisfying
the initial conditions
$$\varphi_0=0,\quad \varphi_1=1,\tag 2.14$$
and it is the analog of the regular solution $\varphi(k,x)$ to the
Schr\"odinger equation
(1.2) with the initial conditions
$$\varphi(k,0)=0,\quad \varphi'(k,0)=1.\tag 2.15$$
Note that $k$ and $-k$ appear in the same way in (1.2) and (2.15), and hence
the regular solution $\varphi(k,x)$ satisfies $\varphi(-k,x)\equiv \varphi(k,x).$
Similarly, $z$ and $z^{-1}$ appear in the same way
in (2.12) and (2.14), and hence
$\varphi_n$ remains unchanged if we replace $z$ with $z^{-1}$ in $\varphi_n.$

Another relevant solution to (2.12) and hence also to (1.1) is
the Jost solution $f_n$ satisfying
the asymptotic condition
$$f_n=z^n[1+o(1)],\qquad n\to+\infty.\tag 2.16$$
It is the analog of the Jost solution $f(k,x)$ satisfying (1.2) with the asymptotics
$$f(k,x)=e^{ ik x}[1+o(1)],\qquad x\to+\infty.\tag 2.17$$

Recall [4,5,10,17,18] that (1.2) has also the solution $g(k,x)$ satisfying the
asymptotics
$$g(k,x)=e^{-ik x}[1+o(1)],\qquad x\to+\infty,$$
and $g(k,x)$ is related to the Jost solution as
$$g(k,x)\equiv f(-k,x).$$
In the discrete case, the corresponding
analogous solution to (2.12) or equivalently to (1.1) is denoted by
$g_n$ and it satisfies the asymptotics
$$g_n=z^{-n}[1+o(1)],\qquad n\to+\infty.\tag 2.18$$
 From (2.12), (2.16), and (2.18) it follows that $g_n$ is obtained
  from $f_n$ by replacing $z$ by
$z^{-1}$ in $f_n.$

When $V_n\equiv 0,$ let us use $\overset{\circ}\to f_n$
and $\overset{\circ}\to g_n$ to denote the corresponding
Jost solution $f_n$ and its relative $g_n,$ respectively.
 From (2.13), (2.16), and (2.18) we see that
$$\overset{\circ}\to f_n=z^n,\quad
\overset{\circ}\to g_n=z^{-n},\qquad n\ge 1,\tag 2.19$$
and from (2.13) and (2.19) we observe that the values of
$\overset{\circ}\to f_n$ and $\overset{\circ}\to g_n$ at $n=0$
are obtained as
$$\overset{\circ}\to f_0=1,\quad
  \overset{\circ}\to g_0=1.\tag 2.20$$

One can directly verify that the regular solution $\overset{\circ}\to \varphi_n$ satisfying
(2.13) and the initial conditions (2.14) is given by
$$\overset{\circ}\to \varphi_n=\ds\frac{z^n-z^{-n}}{z-z^{-1}},\qquad n\ge 0.
\tag 2.21$$
Using (2.9) in (2.21), we can write $\overset{\circ}\to \varphi_n$ also as
$$\overset{\circ}\to \varphi_n=\ds\frac{\sin( n \theta)}{\sin\theta}
,\qquad n\ge 0.\tag 2.22$$
Equivalently, using (2.6) and (2.7) in (2.21) we can express
$\overset{\circ}\to \varphi_n$ in terms of the spectral
parameter $\lambda$ as
$$\overset{\circ}\to \varphi_n=\ds\frac{\left(1-\ds\frac{\lambda}{2}+\ds\frac{1}
{2}\,\sqrt{\lambda(\lambda-4)}\right)^n-
\left(1-\ds\frac{\lambda}{2}-\ds\frac{1}{2}
\,\sqrt{\lambda(\lambda-4)}\right)^n}{\sqrt{\lambda(\lambda-4)}},\qquad n\ge 0.\tag 2.23$$

Let us now investigate the Jost solution $f_n$ to (2.12) when the potential $V$ belongs to class
$\Cal A_b$ specified in Definition~1.1. Using $V_n= 0$ for $n\ge b+1$ in
(2.12) we see that
$$f_n=z^n,\qquad n\ge b.\tag 2.24$$
Let us emphasize that (2.24) holds at $n=b$ as well.
It is convenient to define
$$m_n:= z^{-n} f_n,\tag 2.25$$
which is
the analog of the Faddeev function $m(k,x)$ related to (1.2) and
expressed in terms of the
Jost solution $f(k,x)$ as
$$m(k,x):=e^{-ikx} f(k,x).$$
Note that from (2.24) and (2.25) we get
$$m_n=1,\qquad n\ge b.\tag 2.26$$

Using (2.24) and (2.25) in (2.12) we obtain
$$m_{n-1}=-z^2 m_{n+1}+(z^2+V_n z+1)\,m_n, \qquad n\ge 1,\tag 2.27$$
and hence, we see that (2.27) allows us to obtain the values of
$m_n$ and $f_n$ at $n=0.$

\noindent {\bf Proposition 2.1} {\it Assume that the potential
$V$ belongs to class $\Cal A_b$ specified in Definition~1.1
and that $V_b\ne 0.$ Then:}

\item{(a)} {\it For $0\le n\le b-1,$ the solution $m_n$ to (2.27) with the asymptotic
 condition (2.26) is a polynomial in $z$ of degree $(2b-2n-1)$ with coefficients
 uniquely determined by the ordered
 set $\{V_{n+1},V_{n+2},\dots,V_b\}$ of potential values as}
 $$m_n=\sum_{j=0}^{2b-2n-1} K_{n(n+j)} z^j,\qquad 0\le n\le b-1.\tag 2.28$$
{\it In particular,
for $0\le n\le b-1$ we have}
$$K_{nn}=1,\quad K_{n(n+1)}=\sum_{j=n+1}^b V_j,
\quad K_{n(n+2)}=\ds\sum
_{n+1\le j<l\le b} V_j\, V_l,\tag 2.29$$
$$K_{n(2b-n-2)}=V_b\sum_{j=n+1}^{b-1} V_j, \quad K_{n(2b-n-1)}=V_b
.\tag 2.30$$

\item{(b)} {\it For $0\le n\le b-1,$ the Jost solution $f_n$ to
(2.12) satisfying (2.16) is a polynomial in $z$ of degree $(2b-n-1)$
with the lowest-power term being $z^n$ and
the highest-power term being $V_b z^{2b-n-1},$ and it is given by}
$$f_n=\sum_{j=0}^{2b-2n-1} K_{n(n+j)} z^{n+j}.\tag 2.31$$

\item{(c)} {\it For $0\le n\le b-1,$ the solution $g_n$ to (2.12) satisfying (2.18) is given by}
$$g_n=\sum_{j=0}^{2b-2n-1} K_{n(n+j)} z^{-n-j}.\tag 2.32$$

\item{(d)} {\it The coefficients $K_{nm}$ are real valued, and for
every nonnegative pair of integers $n$ and $m$ we have}
$$K_{nm}=0, \qquad n\ge 2b-m ,\tag 2.33$$
$$K_{nm}=0, \qquad n\ge m+1.\tag 2.34$$

\noindent PROOF: We obtain (2.28) by solving (2.27) iteratively and using (2.26).
We then get (2.31) with the help of (2.25) and (2.28). Finally, (2.32) is obtained by using
$z^{-1}$ instead of $z$ in (2.31). The coefficients $K_{nm}$
appearing in (2.33) and (2.34) are missing in the expansions in
(2.31) and (2.32), and hence they can be chosen as zero.
The real-valuedness of $K_{nm}$ follows from the fact that
$f_n$ and $g_n$ are complex
conjugates of each other when $z$ is on the unit circle $\bold T.$ \qed

Note that, for $1\le n\le b-1$ in going from $f_n$ to $f_{n-1},$ two new powers of $z$ are gained,
which is seen by displaying $f_n$ and $f_{n-1}$ with powers in
$z$ in an ascending order as
$$f_n=z^n+[V_{n+1}+\cdots+V_b]\,z^{n+1}+ \cdots+V_b\,[V_{n+1}+\cdots+V_{b-1}]\, z^{2b-n-2}+V_b\, z^{2b-n-1},$$
$$f_{n-1}=z^{n-1}+[V_n+\cdots+V_b]\,z^n+\cdots+V_b\,[V_n+\cdots+V_{b-1}]\,
z^{2b-n-1}+V_b\, z^{2b-n}.$$
We remark that the coefficients $K_{nm}$ appearing in
Proposition~2.1 will be used in Section~3 in the
Marchenko method for the inverse problem associated with (1.1).

Recall that the Wronskian of any two solutions $\psi(k,x)$
and $\phi(k,x)$ to (1.2) is defined with the help of a matrix determinant as
$$[\psi(k,x);\phi(k,x)]:=\left| \matrix \psi(k,x) & \phi(k,x)\\
\stretch
\psi'(k,x)&\phi'(k,x)\endmatrix \right|.\tag 2.35$$
It is known [10,17,18] that the Wronskian in (2.35) is independent of $x$ and hence
is a function of the spectral parameter $k$ alone.
The analog of (2.35) for the
Wronskian of any two solutions $\psi_n$ and $\phi_n$ to
the discrete equation (1.1) is given by
$$[\psi_n;\phi_n]:=\left| \matrix \psi_n & \phi_n\\
\stretch
\psi_{n+1}&\phi_{n+1}\endmatrix \right|.\tag 2.36$$
Using (2.12) in (2.36) we get
$$\aligned [\psi_n;\phi_n]=&\left| \matrix \psi_n & \phi_n\\
\stretch
-\psi_{n-1}+(z+z^{-1}+V_n) \,\psi_n&-\phi_{n-1}+(z+z^{-1}+V_n) \,\phi_n\endmatrix \right|\\
\stretch
=& \left| \matrix \psi_n&\phi_n \\
\stretch
-\psi_{n-1}  & -\phi_{n-1}\endmatrix \right|.\endaligned\tag 2.37$$
 From (2.37), by interchanging the two
rows in the last determinant and using (2.36), we obtain
$$[\psi_n;\phi_n]=[\psi_{n-1};\phi_{n-1}],$$
confirming that the Wronskian of any two solutions to (1.1)
is independent of $n$ and hence it is
a function of the spectral parameter alone.
Thus, the value of the Wronskian given in (2.36) can be evaluated
at any $n$-value, e.g. at $n=0$ or as $n\to +\infty.$

For $z\ne \pm 1,$ the regular solution $\varphi_n$ can be expressed as a linear combination
of $f_n$ and $g_n$ given in (2.31) and (2.32), respectively. Writing
$$\varphi_n=\alpha\, f_n+\beta\,g_n,\tag 2.38$$
where $\alpha$ and $\beta$ are some coefficients depending only
on the spectral parameter $z,$
we can express $\alpha$ and $\beta$ in terms
of Wronskians as
$$\alpha=\ds\frac{[g_n;\varphi_n]}{[g_n;f_n]},\quad
\beta=\ds\frac
{[f_n;\varphi_n]}{[f_n;g_n]}.\tag 2.39$$
Since the Wronskian of any two solutions to (1.1) is independent
of $n,$ we can evaluate the Wronskians appearing in the
numerators in (2.39) at $n=0,$ and hence with the help of (2.14)
and (2.36) we obtain
$$[g_n;\varphi_n]=g_0,\quad [f_n;\varphi_n]=f_0.\tag 2.40$$
Similarly, we can evaluate the Wronskians appearing in the
denominators in (2.39) as $n\to+\infty,$ and hence with the help of (2.16), (2.18), and
(2.36) we get
$$[f_n;g_n]=z^{-1}-z.\tag 2.41$$
Thus, (2.38) can be written as
$$\varphi_n=\ds\frac{1}{z-z^{-1}}\,\left( g_0 f_n-f_0\, g_n\right).\tag 2.42$$
Note that (2.40) is analogous to
$$[f(-k,x);\varphi(k,x)]=f(-k,0),\quad [f(k,x);\varphi(k,x)]=f(k,0),$$
(2.41) is analogous to
$$[f(k,x);f(-k,x)]=-2ik,$$
and (2.42) is analogous to
$$\varphi(k,x)=\ds\frac{1}{2ik}\left[f(-k,0)
\,f(k,x)-f(k,0)\,f(-k,x)\right].\tag 2.43$$
Let us remark that the property
$\varphi(-k,x)\equiv\varphi(k,x)$ readily follows from (2.43).
Similarly, from (2.42) we directly observe that
replacing $z$ with $z^{-1}$ in $\varphi_n$ does not change
the value of $\varphi_n.$

It is known [10,17,18] that $f(k,0)$ is the so-called Jost function for the Schr\"odinger equation (1.2)
with the Dirichlet boundary condition (2.3). Similarly,
we see that $f_0$ corresponds to the Jost function for the
discrete Schr\"odinger equation (1.1) with
the Dirichlet boundary condition (2.2).
The scattering matrix $S(k)$ for (1.2) with (2.3) is defined as
$$S(k):=\ds\frac{f(-k,0)}{f(k,0)},\tag 2.44$$
and similarly, the scattering matrix associated with (1.1) and (2.2) is
defined as
$$S:=\ds\frac{g_0}{f_0}.\tag 2.45$$
Some relevant properties of the Jost function
$f_0$ are given in the following theorem.

\noindent {\bf Theorem 2.2} {\it Assume that
the potential $V$ belongs to class $\Cal A_b$
specified in Definition~1.1 and that
$V_b\ne 0.$ Then, the Jost function $f_0$ appearing in
(2.42) is a polynomial in $z$ of degree $2b-1,$ it is equal to
$m_0$ appearing in (2.28) with $n=0,$ and it is given by}
$$f_0=\sum_{j=0}^{2b-1} K_{0j} z^j,\tag 2.46$$
{\it where the coefficients $K_{0j}$ are uniquely determined by the ordered set
$\{V_1,\dots,V_b\}$ of potential values.
In particular, we have}
$$K_{00}=1,\quad K_{01}=\sum_{j=1}^b V_j,\quad  K_{02}=\sum_{1\le l<j\le b} V_j V_l ,\tag 2.47$$
$$K_{0(2b-2)}=V_b\sum_{j=1}^{b-1} V_j,\quad
K_{0(2b-1)}=V_b.\tag 2.48$$
{\it Thus, $f_0-1$ is a ``plus" function in the sense that it is analytic in
$z$ when $|z|<1,$ continuous when $|z|\le 1,$ and behaves like $O(z)$ as $z\to 0$ in
$|z|\le 1.$ Similarly, the quantity $g_0$ appearing in (2.42) is given by
$$g_0=\sum_{j=0}^{2b-1} K_{0j} z^{-j},\tag 2.49$$
and hence $g_0-1$ is a ``minus" function in the
sense that it is analytic in
$z$ when $|z|>1,$ continuous when $|z|\ge 1,$ and $O(1/z)$ as $z\to \infty$ in
$|z|\ge 1.$}

\noindent PROOF: With the help of (2.25), the results for $f_0$ directly follow from Proposition~2.1.
The expressions for the
coefficients in (2.47) and (2.48) are obtained with the help of (2.29) and (2.30).
The results for $g_0$ are obtained from $f_0$ given in
(2.46) by replacing $z$ with $z^{-1}.$ \qed

We observe from (2.46) and the first equality in (2.47) that the Jost function
$f_0$ does not vanish at $z=0,$ and in fact the value of the Jost
function at $z=0$ is given by $1.$
It is convenient to display $f_0$ in powers of $z$ in an ascending order as
$$f_0=1+[V_1+\cdots+V_b]\,z+\cdots+V_b\,
[V_1+\cdots+V_{b-1}]\,z^{2b-2}+V_b\,z^{2b-1}.\tag 2.50$$

\noindent {\bf Theorem 2.3} {\it Assume that
the potential $V$ belongs to class $\Cal A_b$ specified in Definition~1.1. Then, the
ordered set $\{V_1,V_2,\dots,V_b\}$ of potential values is uniquely determined
by the ordered set $\{K_{01},K_{12},K_{23},\dots,K_{(b-1)b}\},$
where the quantities $K_{nm}$ are the coefficients appearing in (2.31).}

\noindent PROOF: From the second identity in (2.29) we obtain
$$V_n=K_{(n-1)n}-K_{n(n+1)}, \qquad 1\le n\le b,
\tag 2.51$$
where we also use $K_{b(b+1)}=0,$ which follows from (2.33).
Hence, the proof is complete. \qed

The following theorem describes the behavior of the regular
solution $\varphi_n$ as $n\to+\infty$ when $\lambda=0$ and
when $\lambda=4.$ It shows that (1.1) with the
Dirichlet boundary condition (2.2) cannot have
a bound state at $\lambda=0$ or $\lambda=4.$
It also shows that (1.1) and (1.2) have some similarities
for the generic case and also for the exceptional case.

\noindent {\bf Theorem 2.4} {\it Assume that the potential $V$
belongs to class $\Cal A_b$ specified in Definition~1.1. Let $\lambda$ and
$z$ be the spectral parameters appearing in (1.1) and (2.12),
respectively, and let $\varphi_n$ and $f_n$ be the corresponding
regular solution and the Jost solution to (1.1)
appearing in (2.14) and (2.16), respectively.
Let $\overset{\circ}\to f_n$ and
$\overset{\circ}\to\varphi_n$ be the Jost and regular solutions
appearing in (2.19) and (2.21), respectively,
corresponding to $V_n\equiv 0.$
Then:}

\item{(a)} {\it At $\lambda=0,$ and equivalently at $z=1,$ the
regular solution $\varphi_n$ either grows linearly
as $n\to+\infty,$ which corresponds to the
``generic case," or it remains bounded, which
corresponds to the ``exceptional case." Hence,
$\lambda=0$ never corresponds to a bound state for
(1.1) with the Dirichlet boundary condition (2.2).}

\item{(b)} {\it At $\lambda=4,$ and equivalently at $z=-1,$ the
regular solution $\varphi_n$ generically grows linearly
as $n\to+\infty$ or in the exceptional case it remains bounded. Hence,
$\lambda=4$ never corresponds to a bound state for
(1.1) with the Dirichlet boundary condition (2.2).}

\item{(c)} {\it If $f_0$ has a zero at
$z=1,$ then it must be a simple zero. If $f_0$ has a zero
at $z=1,$ then $g_0$ also has a simple zero at $z=1.$}

\item{(d)} {\it If $f_0$ has a zero at
$z=-1,$ then it must be a simple zero. If $f_0$ has a zero
at $z=-1,$ then $g_0$ also has a simple zero at $z=-1.$}

\item{(e)} {\it The Jost function $f_0$ and the quantity $g_0$
cannot vanish simultaneously at any $z$-value with the possible
exception of
$z=\pm 1.$}

\item{(f)} {\it The exceptional case occurs when
$f_0=0$ at $z=1$ or $f_0=0$ at $z=-1.$ We may have
$f_0$ vanishing at $z=1$ but not at $z=-1,$
we may have $f_0$ vanishing at $z=-1$ but not at $z=1,$
and we may have $f_0$ vanishing at both $z=1$ and $z=-1.$}

\noindent PROOF: For $\lambda=0,$ we can use the method of variation of
parameters [6] and express any solution to (1.1) in terms of two
linearly independent solutions to (2.5). In particular, for
$\lambda=0$ we can use the regular solution $\overset{\circ}\to\varphi_n$
and the Jost solution $\overset{\circ}\to f_n$ that are given by
$$\overset{\circ}\to\varphi_n=n,\quad \overset{\circ}\to f_n=1,\qquad n\ge 0.
\tag 2.52$$
The method of variation of parameters yields the
following expression for the regular solution to (1.1):
$$\varphi_n=n\left[1-\ds\sum_{j=1}^ {n-1} V_j\,\varphi_j
\right]+\ds\sum_{j=1} ^{n-1} j\, V_j\,\varphi_j,\qquad
n\ge 2.\tag 2.53$$
Thus, for $n>b,$ we obtain
$$\varphi_n=n\left[1-\ds\sum_{j=1}^b V_j\,\varphi_j
\right]+\ds\sum_{j=1}^b j\, V_j\,\varphi_j,\qquad
n\ge b+1.\tag 2.54$$
 From (2.54) we conclude that $\varphi_n$ at $\lambda=0$ grows
 linearly in $n$ as $n\to +\infty$ unless the exceptional case
 occurs with
 $$\ds\sum_{j=1}^b V_j\,\varphi_j=1.\tag 2.55$$
 If (2.55) holds, the second summation term in (2.54) must be nonzero
 because otherwise we would have $\varphi_n=0$ for $\lambda=0$
 and $n\ge b+1,$ causing $\varphi_n=0$ for $\lambda=0$ and
 for all $n\ge 0,$
 contradicting the second equality in (2.14).
Thus, in either case, i.e. whether (2.55) holds or not, from
(2.54) we conclude that $\varphi_n$ at $\lambda=0$ cannot
be square summable, and hence there cannot be a bound state
at $\lambda=0$ for (1.1) with the boundary condition (2.2).
For $\lambda=4$ or equivalently $z=-1,$ the proof is obtained
similarly, in which case (2.52) is replaced with
$$\overset{\circ}\to\varphi_n=(-1)^{n-1}
n,\quad \overset{\circ}\to f_n=(-1)^n,\qquad n\ge 0,\tag 2.56$$
(2.53) is replaced with
$$\varphi_n=(-1)^{n-1} n\left[1-\ds\sum_{j=1}^ {n-1}(-1)^j\, V_j\,\varphi_j
\right]+(-1)^{n-1}\ds\sum_{j=1} ^{n-1} (-1)^j\, j\, V_j\,\varphi_j,\qquad
n\ge 2,\tag 2.57$$
(2.54) is replaced with
$$\varphi_n=(-1)^{n-1} n\left[1-\ds\sum_{j=1}^b (-1)^j\, V_j\,\varphi_j
\right]+(-1)^{n-1}\ds\sum_{j=1} ^b (-1)^j\, j\, V_j\,\varphi_j,\qquad
n\ge b+1,\tag 2.58$$
and (2.55) is replaced with
$$\ds\sum_{j=1}^b (-1)^j\, V_j\,\varphi_j=1.$$
Thus, the proofs of (a) and (b) are complete.
Let us now prove (c), (d), and (e).
Since $g_0$ is obtained from $f_0$ by replacing
$z$ by $z^{-1},$ it follows that $g_0$ vanishes at $z=1$ if
$f_0$ vanishes there. Similarly, $g_0$ vanishes at $z=-1$ if
$f_0$ vanishes there. However, $f_0$ and $g_0$ cannot vanish
simultaneously at
any other $z$-value because otherwise, as a result of (2.42), the regular solution
$\varphi_n$ would be zero at that $z$-value for all $n\ge 1,$
contradicting $\varphi_1=1$ stated in (2.14). Thus, in order to
complete the proofs of (c), (d),
and (e), we only need to show that a possible zero
of $f_0$ at $z=1$ must be a simple zero and that
a possible zero
of $f_0$ at $z=-1$ must be a simple zero.
By Proposition~2.1(b) we know that $f_n$ is a polynomial in $z$
and this is true even when $V_b=0.$ Hence, we can expand
$f_n$ around $z=1$ as
$$f_n(z)=f_n(1)+(z-1)\,\dot f_n(1)+O\left((z-1)^2\right),
\qquad z\to 1 \text{ in } \bC,\tag 2.59$$
where $f_n(z)$ denotes the value of $f_n$ at $z$
and $\dot f_n(1)$ denotes the derivative of $f_n$
with respect to $z$ evaluated at $z=1.$
Since $g_n$ is obtained by replacing $z$ with $z^{-1}$ in $f_n,$
 from (2.59) we get
$$g_n(z)=f_n(1)+\left(\ds\frac{1}{z}-1\right)\,\dot f_n(1)+O\left(\left(\ds\frac{1}{z}-1\right)^2\right),
\qquad z\to 1 \text{ in } \bC.\tag 2.60$$
Using
$$\ds\frac{1}{z}-1=-(z-1)+O\left((z-1)^2\right),
\qquad z\to 1 \text{ in } \bC,$$
we can write (2.60) as
$$g_n(z)=f_n(1)-(z-1)\,\dot f_n(1)+O\left((z-1)^2\right),
\qquad z\to 1 \text{ in } \bC.\tag 2.61$$
Using (2.36), (2.59), and (2.61), we write
the left-hand side of (2.41) as
$$[f_n(z);g_n(z)]=2(z-1)\,\left| \matrix
\dot f_n(1) & f_n(1)\\
\stretch
\dot f_{n+1}(1) & f_{n+1}(1) \endmatrix\right|+O\left((z-1)^2\right),
\qquad z\to 1 \text{ in } \bC.\tag 2.62$$
Since we have
$$z-z^{-1}=2(z-1)+O\left((z-1)^2\right),
\qquad z\to 1 \text{ in } \bC,\tag 2.63$$
with the help of (2.41) and (2.63) we obtain
$$[f_n(z);g_n(z)]=-2(z-1)+O\left((z-1)^2\right),
\qquad z\to 1 \text{ in } \bC.\tag 2.64$$
Comparing (2.62) and (2.64) we see that
$$\left| \matrix
f_n(1) & \dot f_n(1)\\
\stretch
f_{n+1}(1) & \dot f_{n+1}(1) \endmatrix\right|=1.\tag 2.65$$
 From (2.65) we conclude that $f_n(1)$ and $\dot f_n(1)$ cannot vanish
simultaneously. Thus, if $f_n(1)=0,$ then we must have $\dot f_n(1)\ne 0.$
Therefore, $z=1$ can only be a simple zero of $f_n(z).$ Proceeding
similarly, we can prove that
if $f_n$ has a zero at $z=-1,$ then such a zero must be a simple zero.
Thus, the proof of (c)-(e) is complete.
Finally, let us prove (f). For this, let us first show that
$\varphi_n$ at $z=1$ remains bounded as $n\to+\infty$
if and only if $f_0(1)=0.$ Using
(2.59) and (2.61) with
$n=0$ as well as
using (2.12), (2.13) and
(2.62), we expand the right-hand side
of (2.42) around $z=1$ as $n\to+\infty,$ and we obtain
$$\varphi_n=n\,f_0(1)-\dot f_0(1)+O\left(z-1\right),
\qquad z\to 1 \text{ in } \bC,\tag 2.66$$
where the terms represented with $O(z-1)$
remain bounded as $n\to+\infty.$
 From (2.66) we conclude that
$\varphi_n$ at $z=1$ remains bounded
if and only if $f_0(1)=0,$ proving that
the exceptional case with $z=1$ occurs if and only if
$f_0(1)=0.$ A similar argument proves that
the exceptional case with
$z=-1$ occurs if and only if
$f_0(-1)=0.$ We can show with some
explicit examples that $f_0$ can
vanish at $z=1$ but not at $z=-1,$ and vice versa.
For example, in Example~6.1 in Section~6,
 from (6.1) we see that we have
 $f_0(1)=0$ and $f_0(-1)=2$ if we choose $V_1=-1.$
 In that example, we get $f_0(-1)=0$ and $f_0(1)=2$ if we choose $V_1=1.$
 In Example~6.3, choosing
 $V_1=-\sqrt{2}$ and
 $V_2=1/\sqrt{2},$ from
 the second equality in (6.12) we see that
$f_0(1)=0$ and $f_0(-1)=0.$
\qed

We remark that
the summation terms appearing in (2.57) and (2.58) are
the discrete analogs of the integrals appearing in (5.5) and (5.6)
of [2].
In the next theorem, we summarize the facts relevant to the bound states of
(1.1) with the boundary condition (2.2). Recall that the
bound states correspond to the $\lambda$-values at which (1.1)
has a square-summable solution satisfying the boundary condition (2.2).

\noindent {\bf Theorem 2.5} {\it Assume that the potential $V$
belongs to class $A_b$ specified in Definition~1.1. Let $\lambda$ and
$z$ be the spectral parameters appearing in (1.1) and (2.12),
respectively, and let $f_n,$ $\varphi_n,$ and $f_0$ be the
corresponding Jost solution
appearing in (2.16), the regular solution appearing in (2.42),
and the Jost solution appearing in (2.46). Then:}

\item{(a)} {\it A bound state can only occur when $\lambda\in (-\infty,0)$
or $\lambda\in (4,+\infty).$ Equivalently, a bound state can only occur
when $z\in(-1,0)$ or $z\in (0,1).$}

\item{(b)} {\it The regular solution evaluated at a
bound state, i.e. the value of $\varphi_n$ evaluated
at a $\lambda$-value corresponding to
a bound state must be real.
Similarly, the Jost solution evaluated at a
bound state must be real.}

\item{(c)} {\it At a bound state the Jost
function $f_0$ has a simple zero in $\lambda$ and
in $z.$ At a bound state the value of
the Jost solution at $n=1$ cannot vanish, i.e. $f_1\ne 0$
at a bound state.}

\item{(d)} {\it The number of bound states $N$ must be finite, and we have
$0\le N\le 2b-1,$ where $b$ is the positive number related to the
support of $V$ and appearing in Definition~1.1. In particular,
we have $N=0$
when $V_n\equiv 0.$}

\item{(e)} {\it At a bound state the quantity
$g_0$ appearing in (2.49) cannot
vanish. Thus, the scattering matrix appearing in (2.45) has a simple
pole at each bound state.}

\noindent PROOF: By its definition, at a bound state, (1.1) must
have a square-summable solution and that solution
must satisfy (2.2). Since the regular solution
$\varphi_n$ satisfies (2.2), any bound-state
solution must be linearly dependent on the
regular solution.
Since the corresponding operator
$\Cal L$ is selfadjoint, a bound state can only occur
when $\lambda$ is real valued. Recall that, as seen from (2.6),
the real $\lambda$-values in the interval $(0,4)$
correspond to the $z$-values on $\bold T^+,$
the upper portion of the
unit circle. On the other hand,
the
$z$-value corresponding to a bound state
cannot occur on $\bold T^+,$ because,
as a consequence of (2.16) and (2.18), neither of the two linearly
independent solutions $f_n$ and $g_n$ to (1.1) vanish at those
$z$-values as $n\to+\infty.$ Furthermore, from Theorem~2.4
we know that a bound state cannot occur
at $z=1$ or at $z=-1.$ We also know from
(2.6) that $z=0$ corresponds to $\lambda=\infty$
and hence a bound state cannot occur at $z=0.$
Thus, we have proved (a).
Let us now prove (b). Since the bound states can only occur at some
real $z$-values, from Proposition~2.1(d) we conclude that $f_n$
given in (2.31) and $g_n$ given in (2.32)
are both real valued at a bound state for any $n\ge 0.$
As a result, from (2.42) we conclude that
the regular solution $\varphi_n$ also
takes real values at a bound state. Thus, we have proved (b).
Let us now turn to the proof of (c).
At a bound state, we conclude
that the coefficient of $g_n$ on the right-hand side
of (2.42) must be zero because otherwise (2.18) implies
that $\varphi_n$ cannot be square summable
at the bound state occurring at some $z$-value
lying in the interval
$z\in(-1,0)$ or $z\in(0,1).$ Thus,
at a bound state the Jost function $f_0$ must be zero.
At a bound state we cannot have $f_1=0$ because
then the only solution to (1.1) with $f_0=0$ and $f_1=0$
would have to be $f_n\equiv 0,$ which
is not compatible with (2.16). Let us now prove that
the zero of $f_0$ at a bound state must be simple. From (1.1)
we obtain
$$\cases
f_{n+1}+f_{n-1}=(2-\lambda+V_n)\, f_n,\\
\stretch
\ds\frac{d f_{n+1}}{d\lambda}+\ds\frac{d f_{n-1}}{d\lambda}=-f_n+
(2-\lambda+V_n)\ds\frac{d f_n}{d\lambda},\endcases
\tag 2.67$$
where the second line is obtained by taking the
$\lambda$-derivative of the first line.
Multiplying the first line in (2.67) by
$df_n/d\lambda$ and the second line by $f_n$ and
taking the difference of the
resulting equations we get
$$f_{n+1}\,\ds\frac{d f_n}{d\lambda}-f_n\,
 \ds\frac{d f_{n+1}}{d\lambda}
 +f_{n-1}\,\ds\frac{d f_n}{d\lambda}-f_n\, \ds\frac{d f_{n-1}}{d\lambda}
 =f_n^2.\tag 2.68$$
 Let us evaluate (2.68) at a $\lambda$-value
 corresponding at a bound state and take the summation
 of both sides starting with $n=1.$ Many cancellations
 on the left-hand side in the summation
 yield
$$f_0\,  \ds\frac{d f_1}{d\lambda}-f_1\, \ds\frac{d f_0}{d\lambda}=
\ds\sum_{n=1}^\infty f_n^2.\tag 2.69$$
 We already know that at a bound state
 $f_0=0,$ $f_1\ne 0,$ and the $f_n$-values are real. Thus,
 (2.69) implies that at a bound state we have
 $$\ds\frac{d f_0}{d\lambda}=-\ds\frac{1}{f_1}\, \ds\sum_{n=1}^\infty f_n^2,$$
 with the right-hand-side being nonzero.
 Thus, the zero of $f_0$ as a function of
 $\lambda$ at the bound state must be simple. From
 (2.6) we see that $d\lambda/dz$ at a bound state
 is strictly positive and hence a simple
 zero of $f_0$ in $\lambda$ corresponds to a simple zero of
 $f_0$ in $z.$ Thus,
 the proof of (c) is complete. By Theorem~2.2 we know that
 $f_0$ is a polynomial in $z$ of degree $2b-1.$
 Since the bound states correspond to the zeros of
 $f_0,$ we conclude that the number of bound states cannot
 exceed $2b-1.$ Furthermore, from the first
 equality in (2.20) we know that the Jost function
 $\overset{\circ}\to f_n$ corresponding to the zero potential
 does not have any zeros
 and hence the number of bound states
 corresponding to the zero potential is zero.
 Thus, we have proved (d).
 At a bound state $z$-value we already know that $f_0$ vanishes
 and that $z$-value cannot be equal to
 $1$ or $-1.$ Thus, from (2.42) we see that the vanishing
 of $g_0$ at a bound state would imply $\varphi_n\equiv 0,$
 contradicting the second condition in (2.14).
Then, from (2.45) and the simplicity of the zeros of
$f_0$ at the bound states, we conclude that
the scattering
matrix $S$ must have a simple pole at each bound state.
Thus, the proof of (e) is complete. \qed

In (2.42) the regular solution $\varphi_n$ to (1.1)
with the initial conditions (2.14) is expressed in
terms of the spectral parameter $z.$ It is possible to express $\varphi_n$
in terms of the spectral parameter $\lambda$ appearing in (1.1).
That result is stated in the next theorem and will be useful in the formulation of
the Gel'fand-Levitan procedure for (1.1).

\noindent {\bf Theorem 2.6} {\it Assume that
the potential $V$ belongs to class $\Cal A_b$ specified in Definition~1.1. Then, for $n\ge 1$ the regular solution
$\varphi_n$ to (1.1) with the initial conditions (2.14)
is a polynomial in $\lambda$ of degree $(n-1)$ and is given by}
$$\varphi_n=\sum_{j=0}^{n-1} B_{nj} \lambda^j,\tag 2.70$$
{\it where the coefficients $B_{nj}$ are real valued and are
uniquely determined by the ordered set
$\{V_1,V_2,\dots,V_{n-1}\}$ of potential values.
In particular, we have}
$$B_{n(n-1)}=(-1)^{n-1},\quad B_{n(n-2)}=(-1)^{n-2} \left[2(n-1)+\sum_{j=1}^{n-1} V_j
\right],
\tag 2.71$$
$$
B_{n(n-3)}=(-1)^{n-3} \left[(n-2)(2n-3)+2(n-2)\sum_{j=1}^{n-1} V_j
+\sum_{1\le k<j\le n-1} V_j V_k\right].\tag 2.72$$

\noindent PROOF:  As seen from (1.1), $\varphi_n$ satisfies
$$\varphi_{n+1}+\varphi_{n-1}=(2-\lambda+V_n)\,\varphi_n,\qquad n\ge 1.
\tag 2.73$$
Solving (2.73) iteratively and using the initial values
$\varphi_0=0$ and $\varphi_1=1$ as stated in (2.14), we get (2.70) with the coefficients
given in (2.71) and (2.72).  Since the initial values are real valued
and the coefficients in (2.73) for $\varphi_{n-1}$ and $\varphi_n$
are real valued for real $\lambda$-values, the
coefficients $B_{nj}$ appearing
in (2.70) are all real valued. \qed

With the help of Theorem~2.6 we see that, for $n\ge 1,$
 the regular solution
$\overset{\circ}\to \varphi_n$ appearing
in (2.23) is a polynomial in $\lambda$ of degree $n-1.$
The same result can also be obtained directly from (2.23)
by using the expansions
$$\left(
1-\ds\frac{\lambda}{2}\pm\ds\frac{1}{2}\,\sqrt{\lambda(\lambda-4)}
\right)^n=\sum_{j=0}^n \binom{n}{j}
\left(1-\ds\frac{\lambda}{2}\right)^{n-j}\left(
\pm\ds\frac{1}{2}\,\sqrt{\lambda(\lambda-4)}
\right)^j,\qquad n\ge 1,\tag 2.74$$
with $\binom{n}{j}$ denoting the binomial coefficient
$n!/(j! \,(n-j)!).$ With the help of (2.74), from (2.23) we
obtain
$$\overset{\circ}\to \varphi_n=\sum_{s=0}^{\lfloor (n-1)/2\rfloor } \binom{n}{2s+1}
\left(1-\ds\frac{\lambda}{2}\right)^{n-1-2s} \lambda^s \left(\ds\frac{\lambda}{4}-1\right)^s,\qquad n\ge 1
,\tag 2.75$$
with $\lfloor x\rfloor$ denoting the floor function, i.e. the greatest
integer less than or equal to $x.$

Using either (2.75) or Theorem~2.6, we directly obtain
the following corollary.

\noindent {\bf Corollary 2.7} {\it For $n\ge 1,$ the regular solution
to (2.5) with the initial conditions (2.14),
namely $\overset{\circ}\to\varphi_n$ appearing in (2.23),
is a polynomial in $\lambda$ of degree $(n-1)$ and is given by}
$$\overset{\circ}\to\varphi_n=\sum_{j=0}^{n-1}
\overset{\circ}\to B_{nj} \lambda^j,\qquad n\ge 1,\tag 2.76$$
{\it where the coefficients $\overset{\circ}\to B_{nj}$ are real valued
and obtained from the expansion in (2.75). In particular, we have}
$$\overset{\circ}\to B_{n(n-1)}=(-1)^{n-1},\quad
\overset{\circ}\to B_{n(n-2)}=  2(-1)^{n-2}(n-1),\tag 2.77$$
$$\overset{\circ}\to B_{n(n-3)}=(-1)^{n-3} (n-2)(2n-3).\tag 2.78$$

Let us now consider expressing the regular
solution $\varphi_n$ to (1.1) in terms of
the elements in the set
$\{\overset{\circ}\to\varphi_j\}_{j=1}^n.$ This problem arises in the formulation of the
Gel'fand-Levitan method, and it is also closely related to the theory of orthogonal
polynomials [9].

\noindent {\bf Theorem 2.8} {\it Assume that
the potential $V$ belongs to class $\Cal A_b$ specified in Definition~1.1. Let $\overset{\circ}\to\varphi_n$ be
the regular solution to (2.5) satisfying the initial
conditions (2.14). Then, for $n\ge 1$ the regular solution
$\varphi_n$ to (1.1) with the initial conditions (2.14) can be written as a linear
combination of the elements in the set $\{\overset{\circ}\to\varphi_j\}_{j=1}^n,$ and we have}
$$\varphi_n=\overset{\circ}\to\varphi_n+ \sum_{j=1}^{n-1} A_{nj}\,\overset{\circ}\to\varphi_j,
\qquad n\ge 1,\tag 2.79$$
{\it with the coefficients
$A_{nj}$ for $j=1,\dots,n-1$ uniquely determined by the ordered set
$\{V_1,\dots,V_{n-1}\}.$ In particular, we have}
$$A_{n(n-1)}=\sum_{j=1}^{n-1} V_j, \qquad n\ge 2,\tag 2.80$$
 {\it with the convention that}
$$A_{n0}= 0,
\quad A_{(n+1)(n+1)}= 1,\qquad n\ge 0,\tag 2.81$$
$$A_{nm}= 0,\qquad 0\le n<m.$$

\noindent PROOF: From Theorem~2.6 we see that
$\varphi_n$ is a polynomial in $\lambda$ of degree
$n-1.$ From Corollary~2.7 we know that,
for $j\ge 1,$ the quantity $\overset{\circ}\to\varphi_j$
is a polynomial in $\lambda$ of degree $j-1.$
Thus, we can use
the set $\{\overset{\circ}\to\varphi_j\}_{j=1}^n$ as a basis
for polynomials in $\lambda$ of degree $n-1,$ and we express
the polynomial $\varphi_n$ as
$$\varphi_n=A_{nn}\overset{\circ}\to\varphi_n+
\sum_{j=1}^{n-1} A_{nj}\,\overset{\circ}\to\varphi_j,\qquad n\ge 1.
\tag 2.82$$
With the help of (2.70) and (2.76), we can compare the coefficients
of $\lambda^{n-1}$ on both sides of (2.82). From the value of
$B_{n(n-1)}$ in (2.71) and the value of
$\overset{\circ}\to B_{n(n-1)}$ in (2.77) we know that those coefficients
are both equal to $(-1)^{n-1},$ and
hence we must have $A_{nn}=1$ in (2.82). Let us remark that
we can use $A_{n0}$ for $n\ge 1$ as the coefficient
of $\overset{\circ}\to\varphi_0,$ by recalling that
$\overset{\circ}\to\varphi_0=0$ as a result of the
Dirichlet boundary condition given in (2.14).
Thus, the choices in (2.81) are appropriate.
Using (2.70), (2.76),
and an analog of (2.76) but with $(n-1)$ instead of $n$
there, we compare the coefficients of
$\lambda^{n-2}$ on both sides of (2.82). This leads us
to the equation
$$B_{n(n-2)}=A_{nn}\overset{\circ}\to B_{n(n-2)}+A_{n(n-1)}
\overset{\circ}\to B_{(n-1)(n-2)},
\qquad n\ge 2.\tag 2.83$$
Inserting in (2.83)
the value $A_{nn}=1$ and the value of $B_{n(n-2)}$ from (2.71)
and the values of $\overset{\circ}\to B_{n(n-2)}$ and $\overset{\circ}\to B_{(n-1)(n-2)}$ from
(2.77), we get (2.80). \qed

We remark that the coefficients $A_{nm}$ will appear
in Section~4 as the unknown quantities in the
Gel'fand-Levitan system for (1.1) with the Dirichlet boundary condition (2.2).
 From (2.80) and the first equality in (2.81) we
obtain the following corollary.

\noindent {\bf Corollary 2.9} {\it Let $V$ be a potential in class $\Cal A_b$
specified in Definition~1.1. Then, $V$
is uniquely
determined by the ordered set $\{A_{21},A_{32},\dots,A_{(b+1)b}\}$
of the corresponding coefficients appearing in (2.79) via}
$$V_n=A_{(n+1)n}-A_{n(n-1)},\qquad n\ge 1,\tag 2.84$$
{\it with the understanding that $A_{10}=0.$}

\vskip 8 pt
\noindent {\bf 3. THE MARCHENKO EQUATION}
\vskip 3 pt

In this section we outline the Marchenko method to recover the potential
$V$ from the corresponding scattering data
for (1.1) with the Dirichlet boundary condition (2.2). We also show that the bound-state
information is contained in the scattering matrix
$S$ appearing in (2.45) as a result of the compact support of
$V.$ Thus, the corresponding scattering matrix alone determines the potential
$V$ in class $\Cal A_b$ uniquely
without any need to specify the bound-state information separately.

The Marchenko equation for (1.2) is given by [5,10,17,18]
$$K(x,y)+M(x+y)+\int_x^\infty dt\,K(x,t)\,M(t+y),\qquad
y>x,\tag 3.1$$
in which the unknown function
$K(x,y)$ is related to the Jost solution $f(k,x)$ appearing in (2.17) as
$$K(x,y):=\ds\frac{1}{2\pi}\int_{-\infty}^\infty dk\,[f(k,x)-e^{ikx}]\,e^{-iky},$$
or equivalently as
$$K(x,y):=\ds\frac{1}{2\pi}\int_{-\infty}^\infty dk\,[f(-k,x)-e^{-ikx}]\,e^{iky}.$$
Note that $K(x,y)= 0$ for $x>y$ as a result of the
facts [5,10,17,18] that
if the potential $V$ in (1.2)
is real valued and integrable on $x\in(0,b)$
and vanishes for $x>b,$ then
for each fixed $x\ge 0$ the Jost solution
$f(k,x)$ is analytic in $k$ in the open
upper-half complex plane $\bCp,$ continuous in $k$ in the
closed upper-half complex plane $\bCpb,$ and
$e^{-ikx} f(k,x)=1+O(1/k)$ as $k\to\infty$ in $\bCpb.$
The kernel of the Marchenko equation is expressed in terms of the scattering matrix
$S(k)$ given in (2.44) and the bound-state information as
$$M(y):=\ds\frac{1}{2\pi}\int_{-\infty}^\infty dk\,[1-S(k)]\,e^{iky}
+\sum_{s=1}^N c_s^2\,e^{-\kappa_s y},\tag 3.2$$
where the values $k=i\kappa_s$ for $s=1,\dots,N$
with distinct positive $\kappa_s$ correspond to the zeros of
$f(k,0)$ and the constants
$c_s$ denote the Marchenko bound-state norming constants given by
$$c_s:=\ds\frac{1}{\ds\sqrt{\int_0^\infty dx\,f(i\kappa_s,x)^2 }}.$$
Because the potential $V$ is assumed to be supported
within the finite interval $x\in [0,b],$ it turns out that
the value of the norming constant $c_s$ is uniquely determined by
the scattering matrix $S(k)$ as [3,4]
$$c_s=\sqrt{i\, \text{Res}[S(k),i\kappa_s]},\tag 3.3$$
with $\text{Res}[S(k),i\kappa_s]$ denoting the residue of $S(k)$ at
the pole $k=i\kappa_s.$ We note that the poles of
$S(k)$ in $\bCp$ are all simple and that the number of bound states
denoted by
$N$ is either $0$ or a positive integer.

One can derive (3.1) from (2.43) as follows. From (2.43) we get
$$f(k,0)\,f(-k,x)-f(-k,0)\,f(k,x)=-2ik\,\varphi(k,x),
\qquad k\in\bR,\tag 3.4$$
where $\bR$ denotes the real axis.
Dividing (3.4) by $f(k,0)$ and using (2.44) we obtain
$$f(-k,x)-S(k)\,f(k,x)=-\ds\frac{2ik}{f(k,0)}\,\varphi(k,x),
\qquad k\in\bR.\tag 3.5$$
We can write (3.5) as
$$[f(-k,x)-e^{-ikx}]+[1-S(k)]\,e^{ikx}+[1-S(k)][f(k,x)-e^{ikx}]=H(k,x),\qquad k\in\bR,
\tag 3.6$$
where we have defined
$$H(k,x):=f(k,x)-
e^{-ikx}-\ds\frac{2ik}{f(k,0)}\,\varphi(k,x).\tag 3.7$$
With the help of the Fourier transform of (3.6) with $(2\pi)^{-1}\int_{-\infty}^\infty dk\,e^{iky}$
for $y>x,$
we obtain (3.1). When the potential $V$ in (1.2) is real valued
and compactly supported, the bound-state information can be obtained
[3,4] from the scattering matrix $S$ given in (2.44) because in that case
$S(k)$ has a meromorphic extension from
$k\in\bR$ to
$k\in\bCp$ with simple poles at $k=i\kappa_s$ for $s=1,\dots,N$
and the residues of $S(k)$ at $k=i\kappa_s$ are related to the
 the bound-state norming constants as in (3.3). Note that
the contribution of the Fourier transform of (3.7) from each bound state
at $k=i\kappa_s$ is evaluated with the help of
(3.4) and the fact that $f(i\kappa_s,0)=0,$ and that
contribution is given by
$$-\ds\frac{1}{2\pi}\int_{-\infty}^\infty
dk\, \ds\frac{2ik}{f(k,0)}\,\varphi(k,x)\,e^{iky}=-i\,\ds\frac{f(-i\kappa_s,0)}
{\dot f(i\kappa_s,0)}\,f(i\kappa_s,x)\,e^{-\kappa_s y},$$
where an overdot denotes the derivative with respect to $k.$ From (2.44) we see that
$$\text{Res}[S(k),i\kappa_s]=\ds\frac{f(-i\kappa_s,0)}{\dot f(i\kappa_s,0)},
\tag 3.8$$
and hence $c_s$ appearing in (3.2) is given by (3.3).

The analog of (3.1) for (1.1) is derived in a similar way. We can write (2.42) as
$$f_0\, g_n-g_0 f_n=-(z-z^{-1})\,\varphi_n,\qquad z\in\bold T,\tag 3.9$$
where we recall that $\bold T$ denotes the unit circle
$|z|=1$ in the complex $z$-plane.
Dividing (3.9) by the Jost function $f_0$ and using (2.45) we get
$$g_n-S f_n=-\ds\frac{z-z^{-1}}{f_0}\,\varphi_n,\qquad z\in\bold T.\tag 3.10$$
We can write (3.10) as
$$(g_n-z^{-n})+(1-S)\,z^n+(1-S)(f_n-z^n)=H_n,\qquad z\in\bold T,\tag 3.11$$
where we have defined
$$H_n:=f_n-z^{-n}-\ds\frac{z-z^{-1}}{f_0}\,\varphi_n.$$
Let us take the Fourier transform of (3.11) with $(2\pi i)^{-1} \oint dz\,z^{m-1}$ for $m\ge n+1,$
where the integral is along the unit circle $\bold T$
 in the counterclockwise direction traversed
once.
This yields the linear system of Marchenko equations
$$K_{nm}+M_{n+m}+\sum_{j=n+1}^{2b-n-1} K_{nj}M_{j+m}=0,\qquad
0\le n<m,\tag 3.12$$
where we have defined
$$K_{nm}:=\ds\frac{1}{2\pi i}\oint dz\,(g_n-z^{-n})\,z^{m-1},\tag 3.13$$
$$M_n:=\ds\frac{1}{2\pi i}\oint dz\,(1-S)\,z^{n-1}+\sum_{j=s}^N c_s^2 \,z_s^n.\tag 3.14$$
In (3.14), the values
$z=z_s$ for $s=1,\dots,N$ correspond to the zeros of $f_0$ inside the unit circle $\bold T$
and the Marchenko bound-state norming constants $c_s$ are
defined as
$$c_s:=\ds\frac{1}{\sqrt{\ds\sum_{n=1}^\infty
f_n^2\big|_{z=z_s}}},\tag 3.15$$
where we have used the fact that
the value of $f_n$ at a bound state
is real as asserted in Theorem~2.5(b),
and hence we have replaced the absolute square in the denominator in
(3.15) with the ordinary square.
Since the potential $V$ is compactly supported, it turns out that the Marchenko norming
constant $c_s$ defined in (3.15) can be evaluated from the scattering matrix
$S$ appearing in (2.45) via
the residue of $S/z$ at the bound-state pole $z=z_s.$ We have
$$c_s=\sqrt{\text{Res}\left[\ds\frac{S}{z},z_s\right]},\qquad s=1,\dots,N.\tag 3.16$$
We recall that each bound state $z=z_s$ occurs when either
$z_s\in (-1,0)$ or $z_s\in (0,1).$
Let us remark that (3.13) can also be written as
$$K_{nm}=\ds\frac{1}{2\pi i}\oint dz\,(f_n-z^n)\,z^{-m-1},\tag 3.17$$
and that the contribution to the Fourier transform from the
right-hand side of (3.11) at each
zero $z=z_s$ of the Jost function $f_0$ is evaluated
by using
$$\aligned\ds\frac{1}{2\pi i}\oint dz\,H_n\,z^{m-1}=&
-\ds\frac{1}{2\pi i}\oint dz\,\ds\frac{z-z^{-1}}{z\, f_0}\,\varphi_n\, z^m\\
\stretch
=&\left. -\ds\frac{g_0}{z_s\, \dot f_0}\, z_s^m\,f_n\right|_{z=z_s},
\endaligned\tag 3.18$$
where we have used (3.9) and the fact that $f_0$ vanishes
at $z=z_s.$
The overdot in (3.18) indicates the derivative
with respect to $z.$
The value of $c_s$ appearing in (3.14)
is given by (3.16) by noting that from (2.45) and (3.18) we have
$$\text{Res}\left[\ds\frac{S}{z},z_s\right]=\left. \ds\frac{g_0}{z_s\,\dot f_0}\right|_{z=z_s},\tag 3.19$$
which is the analog of (3.8) in the discrete case.

Let us note that the quantities $K_{nm}$ given in (3.13) and (3.17)
coincide with those appearing
in (2.28)-(2.32) and (2.46)-(2.49).
The following theorem states that the ordered
set $\{V_1,\dots,V_b\}$ of potential values
is uniquely determined by the scattering
matrix $S$ appearing in (2.45) as well
as by the Jost function $f_0$ appearing in (2.46).

\noindent {\bf Theorem 3.1} {\it Assume that
the potential $V$ belongs to class $\Cal A_b$ specified in Definition~1.1. Then, the Jost function $f_0$ appearing in
(2.46) uniquely determines the potential. Similarly,
the scattering matrix $S$ given in (2.45) uniquely determines the potential. The recovery
of the potential can be accomplished by solving the Marchenko system (3.12) for $K_{nm}$
and then
obtaining the potential $V$ via (2.51).}

\noindent PROOF: Given the Jost function $f_0$ as a function of
$z,$ we also have $g_0$ given in (2.49) because $g_0$ is obtained
by replacing $z$ with $z^{-1}$ in $f_0.$ Thus, as seen from (2.45) we
have at hand the scattering matrix $S$ as a function of the spectral parameter
$z.$ On the other hand, if we are given the scattering matrix $S,$
using (2.45) and (3.16) we then construct the Marchenko kernel
$M_n$ given in (3.14). Next, we solve the Marchenko system
(3.12) and uniquely obtain $K_{nm}$ for $m\ge (n+1)$ for
$n=0,1,\dots,b.$ Finally, we use (2.51) to obtain
the ordered set
$\{V_1,\dots,V_b\}$
of potential values. \qed

\vskip 8 pt
\noindent {\bf 4. SPECTRAL FUNCTIONS AND THE GEL'FAND-LEVITAN METHOD}
\vskip 3 pt

Let us use $\overset{\circ}\to\rho$ to denote the spectral function [9] associated with
the unperturbed discrete Schr\"odinger operator corresponding to
(2.5) with the Dirichlet
boundary condition (2.2). Similarly,
let us use $\rho$ to denote the spectral function associated with
the perturbed discrete Schr\"odinger operator (1.1) with the Dirichlet
boundary condition (2.2).
Our goal in this section is to derive the
corresponding Gel'fand-Levitan system and establish
the recovery of the potential
$V$ in class $\Cal A_b$ from the spectral function $\rho.$
Recall that as the spectral parameter we can use $\lambda$ appearing in (1.1), or
use $z$ appearing in (2.6), or even use $\theta$ appearing in (2.9).
Thus, $\overset{\circ}\to\rho$ and $\rho$ can be viewed as functions of $\lambda$
or functions of $z$ or functions of $\theta.$ We will refer to
$d\overset{\circ}\to\rho$ and $d\rho$ as the spectral measures associated
with (2.5) and (1.1), respectively, with the Dirichlet boundary condition (2.2).

As seen from (1.1)
and (1.3), we observe that the spectral parameter used in [9] is not the same as
the spectral parameter $\lambda$ used in (1.1) but the two are closely related to each other.
Furthermore, the discrete Schr\"odinger equation used in [9] is not the same as
(1.1) we use. Nevertheless, proceeding as in [9], we obtain
our spectral measure $d\overset{\circ}\to\rho$ in terms
of $\lambda$ as
$$d\overset{\circ}\to\rho=\cases 0,\qquad \lambda <0,\\
\stretch
\ds\frac{1}{2\pi}\,\sqrt{\lambda(4-\lambda)}\,d\lambda,\qquad 0\le \lambda\le 4,\\
\stretch
0,\qquad \lambda> 4.\endcases\tag 4.1$$
As a function of $\theta,$ the nonzero portion of
$d\overset{\circ}\to\rho$ corresponding to $0\le\lambda\le 4$
can be expressed as
$$d\overset{\circ}\to\rho=
\ds\frac{2}{\pi}\,\sin^2\theta \,d\theta,\qquad 0\le\theta\le \pi.\tag 4.2$$
As a function of $z,$ the nonzero portion of
$d\overset{\circ}\to\rho$ corresponding to $0\le\lambda\le 4$ is given by
$$d\overset{\circ}\to\rho=-\ds\frac{1}{2\pi i}\,(z-z^{-1})^2\,\ds\frac{dz}{z},
\qquad z\in \overline{\bold T^+},\tag 4.3$$
where we recall that $\overline{\bold T^+}$ denotes the
closure of the upper portion of the
unit circle $\bold T.$
We already know from Theorem~2.5(a) that (2.5) with the Dirichlet boundary condition
(2.2) does not have any bound states and hence the expressions for
$d\overset{\circ}\to\rho$ do not contain any terms related to bound states.

The idea [9,11,16,18] behind the spectral function $\overset{\circ}\to\rho$ is the following. The
regular solution $\overset{\circ}\to\varphi_n$ appearing in (2.21)-(2.23)
as a function of $n$ forms
the sequence $\{\overset{\circ}\to\varphi_n\}_{n=1}^\infty,$
 where each term is a function of the
spectral parameter $\lambda.$ In fact, by Theorem~2.6
we know that, for $n\ge 1,$ the quantity $\overset{\circ}\to\varphi_n$ is
a polynomial in $\lambda$ of degree $n-1.$ The terms in the sequence
$\{\overset{\circ}\to\varphi_n\}_{n=1}^\infty$ form an orthonormal set with respect to the spectral measure
$d\overset{\circ}\to\rho,$ i.e. we have
$$\int d\overset{\circ}\to\rho\,\left( \overset{\circ}\to \varphi_j \, \overset{\circ}\to \varphi_l\right)=\delta_{jl},\tag 4.4$$
where $\delta_{jl}$ is the Kronecker delta
and we omit the integration limits for notational simplicity
with the understanding that the integral is over $\lambda\in\bR.$

The orthonormality condition given in (4.4)
can be stated in $\lambda,$ in $z,$ or in $\theta.$
For example, in terms of $\theta,$ as seen from (2.22) and (4.2) we can write the left-hand side
of (4.4) as
$$\int d\overset{\circ}\to\rho\,\left( \overset{\circ}\to \varphi_j \,
\overset{\circ}\to \varphi_l\right)=\ds\frac{2}{\pi}\int_0^\pi d\theta\, \sin^2\theta\,\left[
\ds\frac{\sin (j\,\theta)}{\sin\theta}\right]
\left[\ds\frac{\sin (l\,\theta)}{\sin\theta}\right].$$
Since we already know that
$$\int_0^\pi d\theta\,\sin (j\,\theta)\, \sin (l\,
\theta)=\ds\frac{\pi}{2}\,\delta_{jl},$$
we see that (4.4) holds when it is expressed in $\theta.$

Let us now consider the main idea behind
the spectral function $\rho$ associated with
(1.1) with the Dirichlet boundary condition (2.2). The regular solution $\varphi_n$ to (1.1)
with the initial conditions (2.14) as a function of $n$ forms
the sequence $\{\varphi_n\}_{n=1}^\infty,$ where
each term in the sequence is a function of the spectral parameter $\lambda.$ In fact,
by Theorem~2.6 we know that for $n\ge 1$ the quantity
$\varphi_n$ is a polynomial in $\lambda$
of degree $n-1.$ The terms in the sequence $\{\varphi_n\}_{n=1}^\infty$
are orthonormal with respect to the spectral measure $d\rho,$ i.e. we have
$$\int d\rho\,\left( \varphi_j \,  \varphi_l\right)=
\delta_{jl},\tag 4.5$$
where again for simplicity we suppress
the integration limits in our notation.

The spectral function $\rho$
contains all the relevant information about the potential $V$
and the boundary condition (2.2).
In the presence of the potential $V,$ assuming that the bound states occur
at $\lambda=\mu_s$ for $s=1,\dots,N$ we have
$$d\rho=\cases \ds\sum_{s=1}^N C_s^2\,\delta(\lambda-\mu_s)\,
d\lambda,\qquad \lambda <0 \text{ or }\lambda>4,\\
\stretch
\ds\frac{1}{2\pi}\,\sqrt{\lambda(4-\lambda)}\,\ds\frac{d\lambda}{
|f_0|^2},\qquad 0\le \lambda\le 4,
\endcases\tag 4.6$$
where $f_0$ is the Jost function appearing in (2.46).
Notice that the continuous part of
$d\rho$ in (4.6) is obtained from $d\overset{\circ}\to\rho$
given in (4.1) via a division by $|f_0|^2.$

Recall
that, as stated in Theorem~2.5(c),
the bound states correspond to the zeros of the Jost function
$f_0.$ The constant $C_s$ appearing
in (4.6) is the Gel'fand-Levitan bound-state norming constant at the bound state
$\lambda=\mu_s,$ i.e.
$$C_s:=\ds\frac{1}{\sqrt{\ds\sum_{n=1}^\infty
\left. \varphi_n^2\right|_{\lambda=\mu_s}
}},\tag 4.7$$
where we have used the fact that
the value of $\varphi_n$ at a bound state
is real, as stated Theorem~2.5(b).

The contribution to the spectral function $\rho$
 from each eigenvalue is given by
 the square of the norming
 constant $C_s$ appearing
 in (4.7). We can obtain the
 spectral density given in (4.6),
 and in particular its continuous part, as
 the limiting case $p\to+\infty$ from the
 spectral density $\rho_p$
 corresponding to the discrete
 Schr\"odinger equation on the finite lattice
 $n\in \{1,2,\dots,p\}$
 with the boundary condition (2.2) at $n=0$
 and an additional boundary condition at $n=p+1,$
 e.g. the Dirichlet condition $\psi_{p+1}=0.$
 The discrete Schr\"odinger operator
 corresponding to such a problem
 on the finite lattice has only the discrete spectrum
 and we can explicitly evaluate the
 corresponding spectral function $\rho_p.$
We then get the continuous part of the
expression in (4.6)
 by taking the limit of
 $\rho_p$ as
 $p\to +\infty.$

By comparing the Marchenko norming constant $c_s$ defined in (3.15)
with the Gel'fand-Levitan norming constant $C_s$ defined in (4.7), we observe that
they differ from each other in the sense that
the Marchenko norming constant $c_s$ is obtained by normalizing the
Jost solution $f_n$ at the bound state $\lambda=\mu_s$ whereas
the Gel'fand-Levitan norming constant is obtained by normalizing
the regular solution $\varphi_n$ at the same bound state. On the other hand,
we already know that the regular solution $\varphi_n$ and the Jost solution
$f_n$ are linearly dependent at a bound state. Thus,
evaluating (2.42) at the bound state $\lambda=\mu_s$ and using the fact that
$f_0$ vanishes at $\lambda=\mu_s$ we obtain
$$\varphi_n\big|_{\lambda=\mu_s}=\left(\ds\frac{g_0}{z-z^{-1}}\bigg|_{\lambda=\mu_s}\right)\,
f_n\big|_{\lambda=\mu_s}.\tag 4.8$$
Since the quantity inside the parentheses in (4.8) is independent of $n,$
by squaring both sides of (4.8) and then using a summation over $n,$ we obtain
$$\ds\sum_{n=1}^\infty \varphi_n^2 \big|_{\lambda=\mu_s}= \left(\ds\frac{g_0}{z-z^{-1}}\right)^2\bigg|_{\lambda=\mu_s}\,\ds\sum_{n=1}^\infty
f_n^2 \big|_{\lambda=\mu_s}.\tag 4.9$$
Using (3.15) and (4.7) in (4.9) we conclude that
the Marchenko and Gel'fand-Levitan bound-state norming constants
are related to each other as
$$c_s=\left|\ds\frac{g_0}{z-z^{-1}}\right|_{\lambda=\mu_s}\,
C_s.\tag 4.10$$
If we know $|f_0|,$ because $g_0$ is obtained from
$f_0$ by replacing $z$ with $z^{-1},$ we can evaluate the factor in (4.10)
relating the Marchenko norming constant $c_s$ and the Gel'fand-Levitan
norming constant $C_s.$
We remind the reader that $z$ and $z^{-1}$ are related to $\lambda$ as in (2.6) and
(2.7), respectively.

As already mentioned, the continuous spectrum of
the discrete Schr\"odinger operator associated with
(1.1) and (2.2) is $\lambda\in [0,4]$ and the discrete
spectrum consists of at most a finite number points $\lambda=\mu_s$
where each $\mu_s$ is located either in the interval
$\lambda\in(-\infty,0)$ or $\lambda\in(4,+\infty).$
Viewed as a function of the spectral parameter
$\theta,$ the continuous part of the spectral
measure $d\rho$ is given by
$$d\rho=\ds\frac{2}{\pi}\,\sin^2\theta \,
\ds\frac{d\theta}{|f_0|^2},\qquad 0\le\theta\le\pi,\tag 4.11$$
which is obtained from (4.6) with the help of the
transformation $\theta\mapsto\lambda$ given
in (2.11).
Similarly, we obtain the same quantity as a function
of the spectral parameter $z$ by using (2.8) in (4.6)
and obtain
$$d\rho= -\ds\frac{1}{2\pi i}\,(z-z^{-1})^2\,\ds\frac{dz}{z\, |f_0|^2},
\qquad z\in \overline{\bold T^+},\tag 4.12$$
Note that, in (4.6), (4.11), and (4.12),
instead of $|f_0|^2,$ we can equally use the product $f_0  g_0$ because
$|f_0|^2=f_0 g_0$ when $z$ is on the upper unit semicircle $\bold T^+$
because $g_0$ is obtained from $f_0$ by
replacing $z$ by $z^{-1},$ which is equivalent to taking the
complex conjugate of $f_0$ when $z$ is on $\bold T^+.$

We derive the linear system of Gel'fand-Levitan equations as follows. Consider
the regular solution $\varphi_n$ to (1.1) with the
initial conditions (2.14).
The set $\{ \varphi_j\}_{j=1}^n$ forms an orthonormal basis
for polynomials in $\lambda$ of degree $n-1$ with respect to the
spectral measure $d\rho$ given in (4.6).
By Theorem~2.6, for $j\ge 1,$ we note that each element $\varphi_j$ in the
set is a polynomial in $\lambda$ of degree $j-1.$
Note that we can write any polynomial in $\lambda$ of degree
$n-2$ or less as a linear combination of the elements
in the set $\{ \varphi_j\}_{j=1}^{n-1}.$
Thus,
with respect to the
spectral measure $d\rho,$ the quantity
$\varphi_n$ is orthogonal to any polynomial with degree $n-2$ or less.
Each element $\overset{\circ}\to\varphi_m$
in the set $\{ \overset{\circ}\to\varphi_m\}_{m=1}^{n-1}$ is a polynomial
in $\lambda$ of degree $m-1,$ and because $m<n$ we can conclude
that $\varphi_n$ is orthogonal to $\overset{\circ}\to\varphi_m$
with respect to the measure
$d\rho,$ i.e. we have
$$\int d\rho\,\left( \varphi_n \, \overset{\circ}\to \varphi_m\right)=0,\qquad
1\le m<n,\tag 4.13$$
where we again suppress the integration limits
for notational simplicity.
Let us replace $d\rho$ in (4.13) with $(d\rho-d\overset{\circ}\to\rho)+d\overset{\circ}\to\rho$
and let us replace $\varphi_n$ in (4.13) with the right-hand side of (2.79). For
$1\le m<n$ we get
$$\aligned 0=& \int (d\rho-d\overset{\circ}\to\rho)\, \overset{\circ}\to \varphi_n \,
\overset{\circ}\to \varphi_m+ \int d\overset{\circ}\to\rho\,
\left( \overset{\circ}\to \varphi_n \,
\overset{\circ}\to \varphi_m\right)\\
\stretch
&
+\sum_{j=1}^{n-1} \int (d\rho-d\overset{\circ}\to\rho)\,A_{nj}\, \overset{\circ}\to \varphi_j \,
\overset{\circ}\to \varphi_m
+\sum_{j=1}^{n-1} \int d\overset{\circ}\to\rho \,
\left( A_{nj}\, \overset{\circ}\to \varphi_j \,
\overset{\circ}\to \varphi_m\right).
\endaligned\tag 4.14$$
Let us define
$$G_{nm}:=\int (d\rho-d\overset{\circ}\to\rho)\, \overset{\circ}\to \varphi_n \,
\overset{\circ}\to \varphi_m.\tag 4.15$$
We will apply the orthonormality given in (4.4) to the second and fourth terms
on the right-hand side in (4.14) and use (4.15) in the first and third terms
on the right-hand side in (4.14). This yields
$$G_{nm}+0+\sum_{j=1}^{n-1} A_{nj}\,G_{jm}+\sum_{j=1}^{n-1} A_{nj}\,\delta_{jm}=0,
\qquad 1\le m< n.\tag 4.16$$
Simplifying (4.16) we obtain the linear system of Gel'fand-Levitan equations
$$A_{nm}+G_{nm}+\sum_{j=1}^{n-1} A_{nj}\,G_{jm}=0,
\qquad 1\le m< n.\tag 4.17$$

The solution to the
inverse problem of recovery of the potential $V$ via the Gel'fand-Levitan method
is achieved as follows. Given the absolute value of the Jost function $f_0$ and
the bound-state data, i.e. the zeros of the Jost function $f_0$ and the corresponding
Gel'fand-Levitan bound-state norming constants $C_s$ appearing in (4.7),
with the help of (4.1) and (4.6) we
form the spectral measures $d\overset{\circ}\to\rho$ and
$d\rho.$ We then evaluate
$G_{nm}$ given in (4.15). Next, we solve the Gel'fand-Levitan system
(4.17) and recover $A_{nm}.$ Then, we use the values of $A_{n(n-1)}$
for $n=2,3,\dots,b$ and recover the ordered set $\{V_1,V_2,\dots,V_b\}$
of potential values via (2.84).

Note that the derivation of the discrete Gel'fand-Levitan system
(4.17) is similar to the following derivation of the continuous
version of the Gel'fand-Levitan equation for (1.2)
with the Dirichlet boundary condition (2.3).
Let $\varphi(k,x)$ be the regular solution
to (1.2) satisfying (2.15). Let us use
$\overset{\circ}\to\varphi(k,x)$ to denote
the regular solution to (1.2) with the boundary condition (2.15)
when $V(x)\equiv 0.$ In fact, we have
$$\overset{\circ}\to\varphi(k,x)=\ds\frac{\sin (kx)}{k},\qquad x\ge 0,$$
where we recall that $\lambda$ and $k$ are related to each other
as $k=\sqrt{\lambda}.$ In the continuous case,
the spectral measures $d\overset{\circ}\to\rho$ and
$d\rho$ are related to the orthonormality
relations given by
$$\int d\overset{\circ}\to\rho\,\left[
\overset{\circ}\to\varphi(k,x)\,
\overset{\circ}\to\varphi(k,y)\right]=\delta(x-y),\qquad
x,y\in[0,+\infty),\tag 4.18$$
$$\int d \rho\,\left[
 \varphi(k,x)\,
 \varphi(k,y)\right]=\delta(x-y),\qquad
x,y\in[0,+\infty),\tag 4.19$$
which are the analogs of (4.4) and (4.5),
respectively. Note that we use
$\delta(x)$ to denote the Dirac delta distribution.
As an analog of (4.15) let us define
$$G(x,y):=\int (d\rho-d\overset{\circ}\to\rho)
\,\left[\overset{\circ}\to\varphi(k,x)\,
\overset{\circ}\to\varphi(k,y)\right].\tag 4.20$$
We have the analog of the expansion in
(2.82), namely
$$\varphi(k,x)=\overset{\circ}\to\varphi(k,x)+
\int_0^x dz\, A(x,z)\, \overset{\circ}\to\varphi(k,z).\tag 4.21$$
 From (4.19) and (4.21) we conclude that, for $y<x,$ the regular solutions
$\varphi(k,x)$ and $\overset{\circ}\to\varphi(k,y)$ are
orthogonal with respect to the spectral measure $d\rho,$ i.e. we have
the analog of (4.13) given by
 $$\int d \rho\,\left[
 \varphi(k,x)\,
 \overset{\circ}\to\varphi(k,y)\right]=0,\qquad
y<x.\tag 4.22$$
Using (4.21) in (4.22) we obtain the analog of (4.14),
namely for $y<x$ we get
$$\aligned 0=& \int (d\rho-d\overset{\circ}\to\rho)
\,
\overset{\circ}\to\varphi(k,x)\, \overset{\circ}\to\varphi(k,y)
+\int d\overset{\circ}\to\rho\,\left[
\overset{\circ}\to\varphi(k,x)\, \overset{\circ}\to\varphi(k,y)\right]\\
& +\int_0^x dz\,\int (d\rho-d\overset{\circ}\to\rho)
\,A(x,z)\,
\overset{\circ}\to\varphi(k,z)\, \overset{\circ}\to\varphi(k,y)+
\int_0^x dz\,\int d\overset{\circ}\to\rho
\,A(x,z)\,
\overset{\circ}\to\varphi(k,z)\, \overset{\circ}\to\varphi(k,y).
\endaligned\tag 4.23$$
Since we use $x>y,$ the second of the four
terms on the right-hand side in (4.23)
vanishes as a result of (4.18). The fourth term there, with the help of
(4.18) is seen to be equal to $A(x,y).$ Furthermore, using (4.20) in the first and
third terms, from (4.23) we get
the Gel'fand-Levitan equation [5,10,17,18]
$$A(x,y)+G(x,y)+\int_0^x dz\, A(x,z)\,G(z,y)=0,
\qquad x>y,\tag 4.24$$
which is the analog of (4.17).
The potential $V$ is recovered from the solution to
(4.24) as
$$V(x)=2\,\ds\frac{d A(x,x)}{dx},$$
which is the analog of (2.84).
Let us finally remark that the analog of
(4.6) in the continuous case is [5,10,17,18]
$$d\rho=\cases
\ds\sum_{s=1}^N C_s^2\, \delta(\lambda-\mu_s),\qquad \lambda <0,\\
\stretch
\ds\frac{1}{\pi}\,
\ds\frac{\sqrt{\lambda}}{|f(\sqrt{\lambda},0)|^2}\, d\lambda,\qquad \lambda\ge 0.
\endcases\tag 4.25$$
where $f(\sqrt{\lambda},0)$ with $k:=\sqrt{\lambda}$
is the Jost function $f(k,0)$
and is obtained from the Jost solution $f(k,x)$ appearing
in (2.17) by putting $x=0.$ In (4.25) we assume
that there are $N$ bound states occurring at
$\lambda=\mu_s$ for $s=1,\dots,N.$ The analog of (4.1)
in the continuous case with
the zero potential is then given by [5,10,17,18]
$$d\overset{\circ}\to\rho=\cases
0,\qquad \lambda <0,\\
\stretch
\ds\frac{\sqrt{\lambda}}{\pi}\, d\lambda,\qquad \lambda\ge 0.
\endcases$$

\vskip 8 pt
\noindent {\bf 5. TRANSMISSION EIGENVALUES}
\vskip 3 pt

In this section we introduce the transmission eigenvalues for
the discrete Schr\"odinger equation (1.1) with the
Dirichlet boundary condition (2.2). We relate the
transmission eigenvalues to the zeros of the
key quantity $D$ defined in (5.2) and also equivalently to
the zeros of the quantity $E$ given in (5.7). We express $D$ in terms of
the spectral parameter $z$ and relate it
as in (5.8) to the Jost function $f_0$ and its relative $g_0,$
 which allows us to obtain the relationship
(5.10) relating
$E$ to $f_0$ and $g_0.$
We introduce the so-called unusual case and characterize it
in several equivalent ways. We recognize that
(5.10) is a discrete analog of a Riemann-Hilbert problem,
which enables us to solve it uniquely in the
usual case and recover $f_0$ from $E.$ Thus, we uniquely
recover the Jost function $f_0$ from the set of transmission eigenvalues
in the usual case.
We are then able to use either of the Marchenko method and the Gel'fand-Levitan method
to uniquely recover the potential
in class $\Cal A_b$ when the transmission eigenvalues
are given as input in the usual case. In the unusual case, the unique recovery may
or may not be possible.
This is illustrated with some examples in Section~6.

The transmission eigenvalues [22] for (1.1) with the
potential $V$ in class
$\Cal A_b$ and with the Dirichlet boundary condition
(2.2) are the $\lambda$-values for which we have nontrivial solutions
$\psi_n$ and $\overset{\circ}\to\psi_n$ to the system
$$\cases -\psi_{n+1}+2 \psi_n-\psi_{n-1}+V_n \psi_n=\lambda \psi_n,\qquad n\ge 1,\\
\stretch
-\overset{\circ}\to\psi_{n+1} +2\overset{\circ}\to\psi_n- \overset{\circ}\to\psi_{n-1}=\lambda
\overset{\circ}\to\psi_n,\qquad n\ge 1,\\
\stretch
\psi_0=0,\quad \overset{\circ}\to\psi_0=0,\\
\stretch
\psi_b=\overset{\circ}\to\psi_b,\quad
\psi_{b+1}=\overset{\circ}\to\psi_{b+1}.
\endcases\tag 5.1$$
In other words, at a transmission eigenvalue $\lambda,$ we must have (1.1)
satisfied for a nontrivial wavefunction $\psi_n,$ the unperturbed problem corresponding to (1.1)
with $V_n\equiv 0$ must have a nontrivial solution
$\overset{\circ}\to\psi_n,$ the solutions $\psi_n$ and $\overset{\circ}\to\psi_n$
both must vanish at $n=0,$ and the solutions $\psi_n$ and $\overset{\circ}\to\psi_n$
must agree at $n=b$ and also at $n=b+1,$ where $b$ is the positive integer
appearing in Definition~1.1 and
related to the support of the potential $V.$ We note that the discrete
transmission-eigenvalue problem described in (5.1) is the analog [3] of the
transmission-eigenvalue problem for the Schr\"odinger equation
with the Dirichlet boundary condition (2.3).

 From the first line of (5.1) and the first equation in the third
line, it follows that $\psi_n$ must be linearly dependent on the
regular solution $\varphi_n$ appearing in (2.42). Similarly,
 from the second line of (5.1) and the second equation in the third
line, it follows that $\overset{\circ}\to\psi_n$ must be linearly dependent on the
regular solution $\overset{\circ}\to\varphi_n$ appearing in (2.23).
The fourth line of (5.1) is equivalent to the statement that the
column vectors $\bm \overset{\circ}\to\psi_b\\
\overset{\circ}\to\psi_{b+1}\endbm$ and $\bm \psi_b\\
\psi_{b+1}\endbm$ are linearly dependent and hence
the matrix formed by using those two vectors
as columns must have zero determinant.
Consequently, the transmission eigenvalues are exactly
the $\lambda$-values for which
the two column vectors $\bm \overset{\circ}\to\varphi_b\\
\overset{\circ}\to\varphi_{b+1}\endbm$ and $\bm \varphi_b\\
\varphi_{b+1}\endbm$ are linearly dependent. Equivalently,
the transmission eigenvalues correspond to the zeros of
the quantity $D$ defined as a matrix determinant via
$$D:=\left| \matrix \overset{\circ}\to\varphi_b &
\varphi_b\\
\stretch
\overset{\circ}\to\varphi_{b+1}&
\varphi_{b+1}
\endmatrix \right|.\tag 5.2$$
When $V_n\equiv 0,$ from (5.1) it follows that
any $\lambda$ is a transmission eigenvalue, and
in this case from (5.2) we get $D\equiv 0.$
When $V_n= 0$ for $n\ge 2$ and $V_1\ne 0,$
 from (5.1) it follows that no
$\lambda$-value can be a transmission eigenvalue, and
in this case from (5.2) we get $D= V_1.$
Thus, in the analysis of the inverse problem
with transmission eigenvalues we can assume that
$b\ge 2,$ where $b$ is the positive integer
related to the support of $V$ given in Definition~1.1.

\noindent {\bf Theorem 5.1} {\it Assume that
the potential $V$ belongs to class $\Cal A_b$ with $V_b\ne 0$ and
further suppose that
$b\ge 2,$ where $b$ denotes the positive integer
related to the support of the potential. Then, the transmission
eigenvalues for (5.1) correspond to the zeros of
the quantity $D$ defined
in (5.2). The quantity
$D$ is a polynomial in $\lambda$ of degree $2b-2,$ and
we have}
$$D=\sum_{s=0}^{2b-2} D_s \lambda^s,\tag 5.3$$
{\it where the coefficients $D_s$
for $s=0,1,\dots,2b-2$ are uniquely determined by
the potential $V.$ In particular, we have}
$$D_{2b-2}=V_b,\quad D_{2b-3}=-V_b\left( 4(b-1)+\sum_{j=1}^{b-1} V_j\right).\tag 5.4$$

\noindent PROOF: The fact that the transmission
eigenvalues correspond to the zeros of
$D$ has already been established in the paragraph
containing (5.2).
Using (2.70)-(2.72) and (2.76)-(2.78) in (5.2), we establish the remaining
results stated. \qed

With the help of Theorem~5.1 we obtain the following result.

\noindent {\bf Theorem 5.2} {\it Assume that
the potential $V$ belongs to class $\Cal A_b$ with $V_b\ne 0$ and
further suppose that
$b\ge 2,$ where $b$ denotes
the positive integer
related to the support of the potential. Then, the quantity $D$
defined in (5.2), as a function of $\lambda,$ has
exactly $(2b-2)$ zeros and hence the number of transmission eigenvalues for
(5.1) is $(2b-2).$ In terms of
the transmission eigenvalues $\lambda_1,\dots,\lambda_{2b-2},$ we have}
$$D=V_b\prod_{j=1}^{2b-2} (\lambda -\lambda_j),\tag 5.5$$
$$\sum_{j=1}^{2b-2}\lambda_j=4(b-1)+\sum_{j=1}^{b-1} V_j.\tag 5.6$$
{\it Hence, the knowledge of the transmission eigenvalues is
equivalent to the knowledge of the quantity $E$ defined as}
$$E:=\ds\frac{D}{V_b},\tag 5.7$$
{\it and the transmission eigenvalues correspond to
the zeros of $E.$}

\noindent PROOF: By Theorem~5.1 we know that the zeros
of $D$ correspond to
the transmission eigenvalues for (5.1) and that
$D$ is a polynomial in $\lambda$ of degree $2(b-1)$
with the leading term $V_b \lambda^{2b-2}.$ Thus, we have established
(5.5). Comparing the coefficients of $\lambda^{2b-3}$ in
(5.3) and (5.5), with the help of (5.4) we get (5.6). \qed

Recall that we can also use $z$ as the spectral parameter
instead of $\lambda,$ where $z$ is related to $\lambda$ as in
(2.6). Next, we express the key quantity $D$ given in (5.2)
in terms of the spectral parameter $z.$ This will help us
to solve the relevant inverse problem to recover the potential
with the help of the Marchenko method or the Gel'fand-Levitan method.

\noindent {\bf Theorem 5.3} {\it Assume that
the potential $V$ belongs to class $\Cal A_b$ with $V_b\ne 0$ and
further suppose that $b\ge 2,$ where $b$ denotes the positive integer
related to the support of the potential. Then:}

\item{(a)} {\it The key quantity $D$ given in (5.2) can be
expressed in terms of the spectral parameter $z$ via}
$$D=\ds\frac{1}{z-z^{-1}}\,(f_0-g_0),\tag 5.8$$
{\it where $f_0$ is the Jost function appearing in (2.46)
and $g_0$ is the quantity in (2.49). The expression (5.8) holds for $b=1$ as well.}

\item{(b)} {\it The quantity $(f_0-1)/V_b$ is uniquely
determined by the transmission
eigenvalues for (5.1). The quantity $K_{01}/V_b$ is also
uniquely determined by
the transmission
eigenvalues, where $K_{01}$ is the coefficient of $z$ in the polynomial
expansion
(2.46) in terms of $z$ for the Jost function $f_0.$ In fact, the
value of $K_{01}/V_b$ is equal to the coefficient of $z$ when
$(f_0-1)/V_b$ is expressed in $z$ as a polynomial.}

\item{(c)} {\it The identity (5.6) is equivalent to}
$$\sum_{j=1}^{2b-2}\lambda_j=4(b-1)+\left(\ds\frac{K_{01}}{V_b}-1\right) V_b.
\tag 5.9$$

\item{(d)} {\it Unless the quantity $K_{01}/V_b$ is equal to $1,$
where the case $K_{01}/V_b=1$ corresponds to the unusual case, $V_b$ is uniquely determined
by the transmission eigenvalues by solving (5.9) for $V_b.$}

\item{(e)} {\it Unless we are in the unusual case,
the transmission eigenvalues uniquely determine
the Jost function $f_0.$ Thus, unless
we are in the unusual case, the transmission eigenvalues uniquely determine
the ordered set  $\{V_1,\dots,V_b\}$ of potential values.}

\noindent PROOF: Using (2.21) and (2.42) in the first and
second columns, respectively, of (5.2) we establish (a).
Using (5.7) in (5.8), we get
$$(z-z^{-1})\, E=\ds\frac{f_0-1}{V_b}-\ds\frac{g_0-1}{V_b}.\tag 5.10$$
 From Theorem~5.2 we know
that the left-hand side of (5.10) is uniquely determined by
the transmission eigenvalues of (5.1). From Theorem~2.2 we know that
$f_0-1$ is a polynomial in $z$ of degree $(2b-1)$ and
the quantity $g_0-1$ is a polynomial in $z^{-1}$ of degree $(2b-1).$
Thus, given the left-hand side of (5.10), we can uniquely
determine $(f_0-1)/V_b$ as follows.
 From Theorem~2.2 it follows that
$(f_0-1)/V_b$ is a ``plus" function and
$(g_0-1)/V_b$ is a ``minus" function, and hence writing
the left-hand side of (5.10) in the spectral
parameter $z,$ we see that the terms containing the
positive powers of $z$ make up
$(f_0-1)/V_b$ and the terms containing the
negative powers of $z$ make up
$(g_0-1)/V_b.$ Having $(f_0-1)/V_b$ at hand, we can use
(3.17) with $n=0$ and $m=1$ to obtain
the value of $K_{01}/V_b.$
We remark that $K_{01}/V_b$ is the coefficient of
$z$ when $(f_0-1)/V_b$ is written in terms of $z$ as a
polynomial in $z.$
Thus, we have proved (b). From the second equality in (2.47) we get
$$\ds\frac{K_{01}}{V_b}=1+\ds\frac{1}{V_b}\sum_{j=1}^{b-1} V_j,\tag 5.11$$
and hence using (5.11) in (5.6) we obtain (5.9), which
establishes (c). By (b) we know that $K_{01}/V_b$ is
uniquely determined by the transmission eigenvalues and
hence we can solve (5.9) for $V_b$ and determine
$V_b$ uniquely in terms of the transmission eigenvalues provided
$K_{01}/V_b\ne 1.$ Thus, we have also established
(d). From (b) we already know that the
transmission eigenvalues uniquely determine
$(f_0-1)/V_b,$ and from (d) we know that in the usual case the
transmission eigenvalues uniquely determine $V_b.$
Thus, in the usual case the
transmission eigenvalues uniquely determine $f_0.$
By Theorem~3.1 we know that $f_0$ uniquely determines the
potential $V.$ Thus, the proof of (e) is complete. \qed

The following result is needed to prove Theorem~5.5. It shows that
a transmission eigenvalue for (5.1) cannot occur at a zero
of the Jost function $f_0.$

\noindent {\bf Proposition 5.4} {\it Assume that
the potential $V$ belongs to class $\Cal A_b$ specified in Definition~1.1.
Let $f_0$ be the Jost function appearing in (2.46),
 $g_0$ the quantity in (2.49), and $D$ the
 quantity defined in (5.2). Let $z$ be the spectral
 parameter appearing in (2.12). Then, the quantity $D$ and $f_0$ cannot vanish simultaneously
at any $z$-value.}

\noindent PROOF: From (5.8) we see that if $D$ and $f_0$ vanished at the same $z$-value with
$z\ne \pm 1,$ then we would also have $g_0$ vanishing at that $z$-value,
which would contradict Theorem~2.4(e). Thus, $D$ and $f_0$ cannot simultaneously
vanish at any $z$-value with $z\ne \pm 1.$ On the other hand, using
(2.59) and (2.61) with $n=0$ as well as (2.41) and
(2.64) in (5.8), we get the expansion for $D$
around $z=1$ as
$$D(z)=\,\dot f_0(1)+O\left(z-1\right),
\qquad z\to 1 \text{ in } \bC,\tag 5.12$$
where we recall that $\dot f_0(1)$
denotes the derivative of
$f_0$ with respect to $z$ evaluated at $z=1.$
Because of (2.65) we know that $\dot f_0(1)\ne 0$ if $f_0(1)=0.$ Thus,
$D$ and $f_0$ cannot simultaneously vanish at $z=1.$ A similar argument
shows that $D$ and $f_0$ cannot simultaneously vanish at $z=-1$ either. \qed

 From (2.20) and (2.45)
 it follows that
 the scattering matrix
 for the unperturbed system
 is given by
$$\overset{\circ}\to S\equiv 1.$$
 The following theorem explores the connection
between the transmission eigenvalues and the energies at which
the scattering from the perturbed system (1.1) and from the
unperturbed system (2.5) coincide. It shows that at a transmission
eigenvalue the scattering matrix $S$ of the perturbed system
takes the value of $1.$ With the possible exception of
$\lambda=0$ and $\lambda=4,$ it shows that any $\lambda$-value
at which $S$ becomes equal to $1$ is a transmission
eigenvalue. It also gives some necessary and sufficient
conditions for each of $\lambda=0$ and $\lambda=4$
to be a transmission eigenvalue for the system (5.1).

\noindent {\bf Theorem 5.5} {\it Assume that
the potential $V$ belongs to class $\Cal A_b$ with $V_b\ne 0$
and further suppose that $b\ge 2,$ where $b$ denotes the positive integer
related to the support of the potential. Then:}

\item{(a)} {\it At a transmission eigenvalue $\lambda$ for (5.1),
the scattering matrix $S$ given in (2.45) takes the value $1.$ Thus,
at a transmission eigenvalue $\lambda$
the scattering from the perturbed system (1.1) with the Dirichlet boundary condition
(2.2) agrees with the
scattering from the unperturbed system (2.5) with the Dirichlet boundary condition
(2.2).}

\item{(b)} {\it With the possible exception of $z=1$ and $z=-1,$
corresponding to $\lambda=0$ and $\lambda=4,$ respectively, each
$\lambda$-value at which the scattering matrix $S$ takes the value of
$1$ is a transmission eigenvalue of (5.1).}

\item{(c)} {\it If the scattering matrix
$S$ takes the value $1$ at $\lambda=0,$ then $\lambda=0$ is a
transmission eigenvalue if and only if
$f_0$ is nonzero and $\dot f_0$ is zero
at $\lambda=0.$}

\item{(d)} {\it If the scattering matrix
$S$ takes the value $1$ at $\lambda=4,$ then $\lambda=4$ is a
transmission eigenvalue if and only if
$f_0$ is nonzero and $\dot f_0$ is zero
at $\lambda=4.$}

\item{(e)} {\it If the scattering matrix
$S$ takes the value $1$ at $\lambda=0,$ then $\lambda=0$ is a
transmission eigenvalue if and only if
$S-1$ has a double zero
at $\lambda=0.$}

\item{(f)} {\it If the scattering matrix
$S$ takes the value $1$ at $\lambda=4,$ then $\lambda=4$ is a
transmission eigenvalue if and only if
$S-1$ has a double zero
at $\lambda=4.$}

\noindent PROOF: Recall that
$\lambda=0,$ as a result of (2.6),
corresponds to $z=1$ in a one-to-one manner and
that $\lambda=4$ corresponds to $z=-1$ in a one-to-one manner.
With the help of (2.45), let us write (5.8) as
$$(z-z^{-1})\,\ds\frac{D}{f_0}=1-S.\tag 5.13$$
By Theorem~5.1 we know that the transmission eigenvalues
correspond to the zeros of $D.$ First,
consider a zero of $D$ occurring at a $z$-value with $z\ne \pm 1.$
 From Proposition~5.4 we know that $f_0$ must be nonzero
at such a $z$-value. Hence, the left-hand side of
(5.13) is zero at that $z$-value, causing $S$ to take the value of
$1$ at that $z$-value. If the zero of $D$ occurs at $z=1,$ then
 from Proposition~5.4 we know that $f_0|_{z=1}\ne 0$ and hence
the left-hand side of (5.13) is zero, which results in $S$
taking the value of $1$ at $z=1.$ A similar argument shows that
if $z=-1$ corresponds to a zero of $D,$ then we have $S$ taking the value of
$1$ at $z=-1.$ Thus, the proof of (a) is complete. The proof of (b)
is given as follows. If $S=1$ at some
$z$-value other than $z=\pm 1,$ then the left-hand side of (5.13) must
be equal to zero at that $z$-value, and hence (5.13) implies that
$D$ must vanish at that same $z$-value. Thus, that $z$-value,
being a zero of $D$ must be a transmission eigenvalue, establishing (b).
To prove (c) we proceed as follows. Using (2.59) and (2.61)
in (2.45), we can expand the scattering matrix $S$ around $z=1$ as
$$S(z)=\ds\frac{f_0(1)-(z-1)\,\dot f_0(1)+O\left((z-1)^2\right) }{  f_0(1)+(z-1)\,\dot f_0(1)+O\left((z-1)^2\right)},
\qquad z\to 1 \text{ in } \bC,$$
which yields, for $z\to 1$ in $\bold C,$
$$S(z)=\cases
1-2 (z-1)\ds\frac{\dot f_0(1)}{f_0(1)}+O\left((z-1)^2\right),\qquad
f_0(1)\ne 0,\\
\stretch
-1+O(z-1),
\qquad f_0(1)=0.\endcases\tag 5.14$$
Thus, we have $S(1)=1$ if $f_0(1)\ne 0$ and we get
$S(1)=-1$ if $f_0(1)= 0.$ On the other hand, from (5.12) we know that
$D$ vanishes at $z=1$ if and only if $\dot f_0(1)=0.$
Therefore, at $z=1$ we have $S=1$ and $D=0$ if and only if $f_0(1)\ne 0$ and $\dot f_0(1)=0.$
Hence, the proof of (c) is complete. The proof of (d) is obtained in a similar way
as in the proof of (c). From (5.14) we see that $S-1$ has a double zero
at $z=1$ if and only if we have
$f_0\ne 0$ and $\dot f_0=0$ at $z=1.$ Thus, (e) is equivalent to (c).
Similarly, one can show that (f) is equivalent to (d), completing the proof
of our theorem. \qed

Let us now investigate the unusual case in the transmission-eigenvalue problem for (5.1). Recall from Theorem~5.3(d) that the unusual case occurs for (5.1) when
$$\sum_{j=1}^{2b-2}\lambda_j=4(b-1),\tag 5.15$$
where $b$ is the integer related to the support of the potential
$V$ and the $\lambda_j$-values for $j=1,\dots,2b-2$
denote the transmission eigenvalues.
 From (5.9) and (5.11),
 or equivalently from (5.6) and (5.15), we see that the unusual case can also be identified
as the case with
$$\sum_{j=1}^{b-1} V_j=0.\tag 5.16$$
In the next theorem we characterize the unusual case in several different ways.

\noindent {\bf Theorem 5.6} {\it Assume that
the potential $V$ belongs to class $\Cal A_b$ with $V_b\ne 0$ and
further suppose that
$b\ge 2,$ where $b$ denotes the positive integer
related to the support of the potential. Then, the unusual case for
the transmission eigenvalue problem for (5.1) occurs if and only if
any one of the following equivalent conditions is satisfied:}

\item{(a)} {\it The transmission eigenvalues $\lambda_j$ satisfies (5.15).}

\item{(b)} {\it The potential $V$ satisfies (5.16).}

\item{(c)} {\it The coefficient $K_{01}$ of $z$ in the polynomial
expression in $z$ given in (2.46) for the Jost function $f_0$
is equal to $V_b.$}

\item{(d)} {\it The coefficient $K_{0(2b-2)}$ of $z^{2b-2}$ in the polynomial
expression in $z$ given in (2.46) for the Jost function $f_0$
is zero.}

\item{(e)} {\it The coefficients of $z$ and $z^{2b-1}$ in the polynomial
expansion in $z$ for $(f_0-1)/V_b$ are both equal to $1.$}

\item{(f)} {\it The coefficient of $z^{2b-2}$ in the polynomial
expansion in $z$ for $(f_0-1)/V_b$
is zero.}

\noindent PROOF: Note that (a) can be used as the definition of the unusual case.
 From (5.9), (5.15), and (5.16) we get the equivalence of (a) and (b).
Using (5.16) in the second equality in (2.47) we see that the unusual
case occurs if and only if $K_{01}=V_b,$ establishing (c). Using (5.16)
in the second equation in (2.48) we observe
that the unusual case is equivalent to having
$K_{0(2b-2)}=0,$ which establishes (d). From (2.46)-(2.48)
it follows that the coefficients of $z$ and $z^{2b-1}$ in the polynomial
expansion in $z$ for $(f_0-1)/V_b$ are given by $K_{01}/V_b$ and $K_{0(2b-1)}/V_b,$
respectively, and hence with the help of (5.11) and (5.16) we see that
$K_{01}/V_b=1$ in the unusual case and $K_{0(2b-1)}/V_b=1$ always, which
establishes (e). From (2.46) and the first equality in (2.48),
 or equivalently from (2.50), we see that the coefficient
of $z^{2b-2}$ in the polynomial
expansion in $z$ for $(f_0-1)/V_b$
is zero if and only if (5.16) holds, establishing (f). \qed

In the usual case we already know that the transmission
eigenvalues uniquely determine the potential.
By presenting some explicit examples, we will show
in the next section that the transmission eigenvalues
for (1.1) may or may not uniquely determine the potential $V$ in the unusual case.
We will show that in the unusual case we may have a one-parameter family
of real-valued potentials corresponding to a given set of transmission eigenvalues,
we may have a finite number of distinct real-valued potentials
corresponding to a given set of transmission eigenvalues, or we may have a
unique real-valued potential corresponding to a given set of transmission eigenvalues.

\vskip 8 pt
\noindent {\bf 6. EXAMPLES}
\vskip 3 pt

In this section we present various explicit examples to illustrate the results
of the previous sections. In particular, we illustrate the recovery
of the potential from the transmission eigenvalues, the
Marchenko method, the Gel'fand-Levitan method, the unusual case, and the bound states.

In the first example below we illustrate the fact that the bound-state information is contained
in the scattering matrix $S$ when the corresponding
potential $V$ belongs to class $\Cal A_b$ described in Definition~1.1.
As stated in Theorem~2.5(e) the bound-state energies
correspond to the poles of the scattering matrix, and
as indicated in (3.19)
the corresponding Marchenko bound-state norming constants are
obtained from $S$ via a residue evaluation.

\noindent {\bf Example 6.1} Assume that the potential is supported at $n=1$
with $V_1\ne 0,$ and hence
$V_n= 0$ for $n\ge 2.$ From (2.24) we have
$f_n=z^n$ for $n\ge 1.$ Using (2.29) and (2.46) we obtain the Jost function as
$$f_0=V_1 z+1.\tag 6.1$$
Thus, the only zero of $f_0$ occurs at $z=-1/V_1.$ From Theorem~2.5(a)
we know that
a bound state exists if the zero of $f_0$ is in
one of the real line segments $z\in(0,1)$ and $z\in (-1,0).$
Thus, there is exactly one bound state if $V_1<-1$ or $V_1>1.$
Assume that one of these conditions
is satisfied so that we have a bound state at $z=-1/V_1.$
 From (6.1) we then get
$$f_n\big|_{z=-1/V_1}=\left(-\ds\frac{1}{V_1}\right)^n.\tag 6.2$$
Using (6.2) in (3.15) we evaluate the Marchenko bound-state norming constant as
$$c_1=\ds\frac{1}{\ds\sqrt {\sum_{n=1}^\infty \ds\frac{1}{V_1^{2n}}}},$$
which simplifies to
$$c_1=\sqrt{V_1^2-1}.\tag 6.3$$
 From (6.1), replacing $z$ by $z^{-1},$ we get the quantity $g_0$ appearing in
(2.49) as
$$g_0=V_1 z^{-1}+1.\tag 6.4$$
Using (6.1) and (6.4) in (2.45), we obtain the scattering matrix $S$ as
$$S=\ds\frac{ V_1 z^{-1}+1}{V_1 z+1}.\tag 6.5$$
Let us verify that (3.16) yields the Marchenko bound-state norming constant
given in (6.3). From (6.5) we evaluate the residue of $S/z$ at
$z=-1/V_1,$ and we have
$$\text{Res}\left[\ds\frac{S}{z},-\ds\frac{1}{V_1}\right]=V_1^2-1,$$
and a comparison with (6.3) shows that (3.16) is verified.
Regarding the transmission eigenvalues, from (5.3) we get
$D\equiv V_1$ and hence there are no transmission eigenvalues.

In the second example given below, we demonstrate that the transmission eigenvalues
may have multiplicities greater than one. We also illustrate that in the unusual case
there may be infinitely many distinct potentials corresponding to
a given set of transmission eigenvalues.

\noindent {\bf Example 6.2} Let us assume that $V_n= 0$ for $n\ge 1$ except
when $n=b$ for some fixed positive integer $b$ and that $V_b$
is a nonzero real parameter. By Theorem~5.6(b) we
see that this corresponds to a special case of the unusual case.
Proceeding as in the proof of Theorem~2.6 we can determine all the values of
the regular solution $\varphi_n.$ We obtain
$$\varphi_b=
\overset{\circ}\to\varphi_b,\quad
\varphi_{b+1}=
\overset{\circ}\to\varphi_{b+1}+V_b \,\overset{\circ}\to\varphi_b,\tag 6.6$$
where we recall that $\overset{\circ}\to\varphi_n$ is the regular solution
when $V_n\equiv 0$ and it appears in (2.21) and (2.22).
Using (6.6) in (5.2), we simplify the determinant defining $D$ and get
$$D=\left|\matrix \overset{\circ}\to\varphi_b & 0\\
\stretch
\overset{\circ}\to\varphi_{b+1}& V_b\, \overset{\circ}\to\varphi_b\endmatrix \right|.
\tag 6.7$$
 From (5.7) and (6.7) we observe that the quantity $E$ is given by
$$E=\left(\overset{\circ}\to\varphi_b\right)^2.\tag 6.8$$
Using (2.22) in (6.8) we have
$$E=\ds\frac{\sin^2 (b\,\theta)}{\sin^2\theta}.\tag 6.9$$
Comparing (5.5) and (6.9) we see that there are $2(b-1)$
transmission eigenvalues $\lambda_j$
including multiplicities, and in fact the transmission eigenvalues
occur when $\theta=j\pi/b$ for $j=1,2,\dots,b-1$
and each has multiplicity two. Using (2.11) we obtain
the $\lambda$-values for the transmission eigenvalues as
$$\lambda_j=2-2\cos\left(\ds\frac{j\pi }{b}\right),\qquad j=1,2,\dots,b-1.\tag 6.10$$
Thus, with the help of (5.5), (5.7), and (6.10) we can write
(6.9) as
$$E=\ds\prod_{j=1}^{b-1} \left[\lambda-2+2\cos\left(\ds\frac{j\pi }{b}\right)\right]^2
\tag 6.11$$
 From (6.11) we observe that there are no transmission eigenvalues if $b=1,$
there is a transmission eigenvalue at $\lambda=2$ with multiplicity two
if $b=2,$ and there are
the transmission eigenvalues at $\lambda=1$ and $\lambda=3$ with multiplicities
two if $b=3.$ It is clear from (6.11) that the value of $V_b$ cannot be determined from the
transmission eigenvalues, and hence there is a one-parameter
family of potentials corresponding to the transmission eigenvalues
for $b\ge 2,$ with $V_b$ being the
nonzero real parameter.

In the third example given below we illustrate that
the transmission eigenvalues may be all real, may be all complex, or may
have multiplicities greater than one.

\noindent {\bf Example 6.3} Assume that $b=2$
in Definition~1.1 and hence $V_n= 0$ for $n\ge 3.$
Further assume that $V_2\ne 0.$
Thus, the potential is known by specifying $V_1$ and $V_2.$
 From (2.24) and (2.29)-(2.31) we get
$f_n=z^n$ for $n\ge 2$ and
$$f_1=V_2\,z^2+z ,\quad f_0= V_2\,z^3 +V_1 V_2 z^2+(V_1+V_2)\,z+1.\tag 6.12$$
Note that we can get $f_0$ given in (6.12) also from (2.50).
Thus, the quantity $g_0$ appearing in (2.49) is obtained from $f_0$ given in
(6.12) by replacing $z$ with $z^{-1}.$ We have
$$g_0= V_2\,z^{-3} +V_1 V_2 z^{-2}+(V_1+V_2)\,z^{-1}+1.\tag 6.13$$
Using (6.12) and (6.13) in (5.8) we obtain
$$D=\ds\frac{ V_2\,(z^3-z^{-3}) +V_1 V_2 (z^2-z^{-2})+(V_1+V_2)\,(z-z^{-1})}{z-z^{-1}},$$
or equivalently
$$D=V_2(z^2+1+z^{-2})+V_1 V_2 (z+z^{-1})+(V_1+V_2).\tag 6.14$$
Using (2.8) in (6.14) we express $D$ as a function of
$\lambda$ as
$$D=V_2 \lambda^2 -V_2 (4+V_1) \lambda+(V_1+2V_1 V_2+4 V_2).\tag 6.15$$
Thus,
the quantity $E$ is obtained by using (6.15) in (5.7) as
$$E=\lambda^2 -(4+V_1) \lambda+\left(4+ 2V_1+\ds\frac{V_1}{V_2}\right).\tag 6.16$$
We know by Theorem~5.2 that the transmission eigenvalues
correspond to the zeros of $E$ given in (6.16). Thus, from (6.16) we see that
the two transmission
eigenvalues are given by
$$\lambda=2+\ds\frac{V_1}{2}\pm\sqrt{\ds\frac{V_1^2}{4}-\ds\frac{V_1}{V_2}}.
\tag 6.17$$
Hence, depending on the potential,
we may have two real distinct transmission eigenvalues,
a double real transmission eigenvalue, or a pair of
complex conjugate transmission eigenvalues. For example,
we have a double transmission eigenvalue at
$\lambda=0$ if $(V_1,V_2)=(-4,-1),$
we have $(\lambda_1,\lambda_2)=(i,-i)$
when $(V_1,V_2)=(-4,-4/5),$
we have $(\lambda_1,\lambda_2)=(3,-3)$ when $(V_1,V_2)=(-4,4/5),$
and we have $(\lambda_1,\lambda_2)=(1+i,1-i)$ when $(V_1,V_2)=(-2,-1).$
 From Theorem~5.6(b) it follows that the unusual case occurs when $V_1=0,$ in which
case we have
$$E=(\lambda-2)^2.$$
As already explained in Example~6.2,
in the unusual case, corresponding to the
transmission eigenvalues $(\lambda_1,\lambda_2)=(2,2)$
we have a one-parameter family
of potentials with $V_1=0$ and $V_2$ being the nonzero real parameter.

In the fourth example given below, we illustrate Theorem~5.5(c) and Theorem~5.5(d).

\noindent {\bf Example 6.4} As in Example~6.3, let us assume
that $b=2$ with $V_2\ne 0,$ and hence $V_n= 0$ for $n\ge 3.$
Let us first illustrate Theorem~5.5(c). Using $f_0$ given in (6.12),
at $z=1$ we obtain
$$f_0(1)=2V_2+V_1+V_1 V_2+1,\quad \dot f_0(1)=4 V_2+V_1+2 V_1 V_2.\tag 6.18$$
 From (6.18) we see that $\dot f_0(1)=0$ is equivalent to
$$V_2=-\ds\frac{V_1}{4+2 V_1},\tag 6.19$$
in which case from (5.14) we get $S|_{z=1}=1.$
Using (6.19) in the first equality in (6.18), we see that $\dot f_0(1)\ne 0$
is equivalent to $V_1\ne -2.$ From the denominator on the right-hand side in
(6.19) we already see that we must have $V_1\ne -2$ for (6.19) to make sense.
Let us now determine when $\lambda=0$ becomes
a transmission eigenvalue. From (6.14) we see
that $\lambda=0$ is a transmission eigenvalue if and only if the sum of the last
three terms on the right-hand side in (6.14) is zero, i.e.
$$4+2 V_1+\ds\frac{V_1}{V_2}=0,$$
which is equivalent to having (6.19). Thus, $\lambda=0$ is a transmission eigenvalue
if and only if
we have $f_0(1)\ne 0$ and $\dot f_0(1)=0,$ verifying Theorem~5.5(c).
Let us now illustrate Theorem~5.5(d). Evaluating $f_0$ of (6.12)
at $z=-1$ we get
$$f_0(-1)=-2V_2-V_1+V_1 V_2+1,\quad \dot f_0(-1)=4 V_2+V_1-2 V_1 V_2.\tag 6.20$$
 From (6.20) we see that $\dot f_0(-1)=0$ is equivalent to
$$V_2=-\ds\frac{V_1}{4-2 V_1}.\tag 6.21$$
Using (6.21) in the first equality in (6.20), we see that $f_0(-1)\ne 0$
is equivalent to $V_1\ne 2.$ From the denominator on the right-hand side in
(6.21) we already see that we must have $V_1\ne 2$ for (6.21) to make sense.
Let us now determine when $\lambda=4$
becomes a transmission eigenvalue. From (6.17) we see
that $\lambda=4$ is a transmission eigenvalue if and only if
$$4=2+\ds\frac{V_1}{2}\pm\sqrt{\ds\frac{V_1^2}{4}-\ds\frac{V_1}{V_2}},$$
or equivalently if and only if
$$\left(-2+\ds\frac{V_1}{2}\right)^2=\ds\frac{V_1^2}{4}-\ds\frac{V_1}{V_2}.
\tag 6.22$$
Simplifying (6.22) we see that (6.22) is equivalent to (6.21).
Thus, $\lambda=4$ is a transmission eigenvalue
if and only if we have
$f_0(-1)\ne 0$ and $\dot f_0(-1)=0,$ verifying Theorem~5.5(d).

In the fifth example given below we illustrate that, in the unusual case,
we may have a unique real-valued potential $V$ in class
$\Cal A_b$ corresponding to a given
set of transmission eigenvalues,
we may have a finite number of distinct real-valued potentials
corresponding to a given set of transmission eigenvalues, or
we may have infinitely many distinct real-valued potentials corresponding to a given set of transmission eigenvalues.

\noindent {\bf Example 6.5} Assume that $b=3$ with $V_3\ne 0,$
and hence $V_n= 0$ for
$n\ge 4.$ By Theorem~5.2, the system
(5.1) has four transmission eigenvalues uniquely determined by the ordered set
$\{V_1,V_2,V_3\}.$ By Theorem~5.6(b) we know that the unusual case
occurs when $V_2=-V_1.$ Let us see to what extent the
transmission eigenvalues determine the potential $V$
in the unusual case.
Proceeding as in the proof of Theorem~2.6, we can evaluate the values of
the regular solution $\varphi_n.$ We get
$$\varphi_3=\lambda^2 -4\lambda+3-V_1^2,\tag 6.23$$
$$\varphi_4=
-\lambda^3+(6+v_3)\lambda^2+(V_1^2-4V_3-10)\lambda+(4-V_1-2V_1^2+3V_3-V_3 V_1^2),
\tag 6.24$$
and hence we also have
the values of the regular solution $\overset{\circ}\to\varphi_n$
corresponding to $V_n\equiv 0,$ and in particular from (2.75) we obtain
$$\overset{\circ}\to\varphi_3=\lambda^2 -4\lambda+3,\quad
\overset{\circ}\to\varphi_4=
-\lambda^3+6\lambda^2-10\lambda+4.
\tag 6.25$$
Since $b=3$ in our example, using (6.23)-(6.25) in (5.2), with the help of
(5.7) we get
$$E=\lambda^4-8\lambda^3+(22-\gamma)\lambda^2+
(4\gamma+\epsilon-24)\lambda+(9-3\gamma-2\epsilon),
\tag 6.26$$
where we have let
$$\gamma:=V_1^2+\ds\frac{V_1}{V_3},\quad \epsilon:=\ds\frac{V_1^2}{V_3}.\tag 6.27$$
Recall that we assume $V_2=-V_1,$ which signifies the unusual case.
Since the inverse problem of determining the potential $V$ from
the transmission eigenvalues is equivalent to
the recovery of $V$ from $E,$ our inverse problem is equivalent to
the determination of $V_1$ and $V_3$ from the ordered set $\{\gamma,\epsilon\}.$
Let us analyze this problem algebraically and see to what extent we can
determine the real parameter $V_1$ and the nonzero real parameter
$V_3$ from the ordered set $\{\gamma,\epsilon\}.$ Eliminating
$V_3$ in (6.27) we obtain
$$V_1^3-\gamma \,V_1+\epsilon=0.\tag 6.28$$
Given the ordered set $\{\gamma,\epsilon\},$ the solution to the algebraic
equation in (6.28) may yield one, two, or three distinct
real values for $V_1.$ Once we obtain $V_1$ by solving
(6.28), we can recover $V_3$ from the second identity in (6.27).
If $\gamma=0$ and $\epsilon=0,$ then $V_1=0$ is the solution with
multiplicity three. This special case is analyzed in Example~6.2
already, and it corresponds to the unusual case.
We already know that in this unusual case, $V_3$ cannot be determined
 from the knowledge of transmission eigenvalues.
In this particular case the four transmission eigenvalues are
given by
$$\lambda_1=1,\quad \lambda_2=1,\quad
\lambda_3=3,\quad
\lambda_4=3,$$
corresponding to the one-parameter family of potentials
with $V_n= 0$ for $n\ne 3$ and
$V_3$ being any nonzero real constant.
Let us next consider the case $\gamma=0$ and $\epsilon\ne 0.$
In this case from (6.27) and (6.28) we get
$$V_1=-\root{3}\of {\epsilon},\quad V_3=\ds\frac{1}{\root{3}\of {\epsilon}},$$
indicating that the transmission eigenvalues uniquely determine
the potential $V$ in this particular case.
For example, the choice $\gamma=0$ and $\epsilon=1$ yields
the unique potential $V$ with
$$(V_1,V_2,V_3)=(-1,1,1),$$
corresponding to the transmission eigenvalues
$$\lambda_1=1.4751\overline{1},\quad \lambda_2=0.50978\overline{4},\quad
\lambda_3=3.0075\overline{5}+0.51311\overline{6}\, i,\quad
\lambda_4=3.0075\overline{5}-0.51311\overline{6}\, i,$$
which are obtained as the zeros of $E$ given in (6.26).
Note that we use an overline on a digit to indicate a round off.
Let us next consider the case where
(6.28) has three distinct real solutions for $V_1.$
For example, when
$\gamma=7$ and $\epsilon=6,$ we have three distinct solutions for
$V_1$ in (6.28), namely, $1,$ $2,$ and $-3.$ In this case, the transmission
eigenvalues $\lambda_j$ are obtained as the zeros of $E$ from (6.26) as
$$\lambda_1=4,\quad \lambda_2=3.8557\overline{7},\quad
\lambda_3=1.3216\overline{4},\quad
\lambda_4=-1.1774\overline{1},$$
corresponding to exactly the three distinct potentials with the values of
$(V_1,V_2,V_3)$ given by
$$\left(1,-1,\ds\frac{1}{6}\right), \quad \left(2,-2,
\ds\frac{2}{3}\right),\quad \left(-3,3,\ds\frac{3}{2}\right).$$
When $\gamma=3$ and $\epsilon=2,$ we get exactly two distinct real
$V_1$-values as solutions to (6.28), namely $V_1=1$ and $V_1=-2.$
In this particular case, the transmission
eigenvalues $\lambda_j$ are obtained as the zeros of $E$ from (6.26) as
$$\lambda_1=1.249\overline{5},\quad \lambda_2=-0.25862\overline{7},\quad
\lambda_3=3.5045\overline{6}+0.30989\overline{2}\, i,\quad
\lambda_4=3.5045\overline{6}-0.30989\overline{2}\, i,$$
corresponding to exactly the two distinct potentials with the values for
$(V_1,V_2,V_3)$ given by
$$\left(1,-1,\ds\frac{1}{2}\right), \quad \left(-2,2,2\right).$$
As a summary, in this example, corresponding to the four transmission
eigenvalues, we have demonstrated that
we may have a unique potential in class $\Cal A_b,$ two distinct
potentials, three distinct potentials, or a
one-parameter family of infinitely many potentials.

In the sixth example below we illustrate the use of
the Marchenko method to obtain the potential from the transmission
eigenvalues.

\noindent {\bf Example 6.6} Assume that $b=2$ with $V_2\ne 0,$
and hence $V_n= 0$ for
$n\ge 3.$ We have $2(b-1)$ equal to $2,$ and hence
there are two transmission eigenvalues, as indicated in Theorem~5.2. Let us assume that they are
given by
$$\lambda_1=\ds\frac{11+\sqrt{57}}{4},\quad \lambda_2=\ds\frac{11-\sqrt{57}}{4}.
\tag 6.29$$
Using (5.5), (5.7), and (6.29) we see that
the quantity $E$ defined in (5.7) is given by
$$E=\left(\lambda-\ds\frac{11+\sqrt{57}}{11}\right)
\left(\lambda-\ds\frac{11-\sqrt{57}}{11}\right).
\tag 6.30$$
Using (2.8) in (6.30) we express $E$ in terms of
the spectral parameter $z$ and obtain
$$E=z^2+\ds\frac{3}{2} \,z-1 +\ds\frac{3}{2}\, z^{-1}+z^{-2}.$$
Thus, we have
$$(z-z^{-1})\,E=z^3+\ds\frac{3}{2} \,z^2-2z+2z^{-1}-\ds\frac{3}{2} \,z^{-2}+z^{-3}.\tag 6.31$$
 From (6.31), with the help of (5.10)
we extract $(f_0-1)/V_2$ by taking the positive powers of $z$ and get
$$\ds\frac{f_0-1}{V_2}=  z^3+\ds\frac{3}{2} \,z^2-2z.\tag 6.32$$
As seen from Theorem~5.3(b), the coefficient of $z$ on the right-hand
side in (6.32) gives us $K_{01}/V_2,$ and hence
$$\ds\frac{K_{01}}{V_2}=-2.\tag 6.33$$
Using (6.29) and (6.33) in (5.9), we write (5.9) as
$$\ds\frac{11}{2}=4+\left(-2-1\right)\,V_2,$$
and hence recover $V_2$ as
$$V_2=-\ds\frac{1}{2}.\tag 6.34$$
Note that the first equality in (2.48) indicates that
the coefficient of $z^2$ in (6.32) must be $V_1,$ i.e. we have
$$V_1=\ds\frac{3}{2}.$$
However, if $b$ were a large integer, we would not be able to obtain
the ordered set $\{V_1,V_2,\dots,V_{b-1}\}$ of potential values
in such a straightforward manner. Thus,
we will use the Marchenko procedure to
illustrate the recovery of $V_1$ because
it is a more systematic way of recovering the potential.
Using (6.34) in (6.32) we obtain $f_0$ as
$$f_0=-\ds\frac{1}{2}\, z^3-\ds\frac{3}{4}\, z^2+z+1.\tag 6.35$$
 From (2.46) and (2.49) we see that the quantity $g_0$ is obtained
by replacing $z$ with $z^{-1}$ in
(6.35) and we have
$$g_0=-\ds\frac{1}{2}\, z^{-3}-\ds\frac{3}{4}\, z^{-2}+z^{-1}+1.\tag 6.36$$
Using (6.35) and (6.36) in (2.45) we obtain the scattering matrix $S$ as
$$S=\ds\frac{-\ds\frac{1}{2}\, z^{-3}-\ds\frac{3}{4}\, z^{-2}+z^{-1}+1}{
-\ds\frac{1}{2}\, z^3-\ds\frac{3}{4}\, z^2+z+1}.\tag 6.37$$
Let us illustrate the determination of the potential $V$ from $S$
given in (6.37).
We determine the zeros of $f_0$ inside the unit circle $|z|=1,$
which is equivalent to finding the poles of $S,$
as indicated in Theorem~2.5(e). Note
that $f_0$ has three zeros at $z=z_s$ for $s=1,2,3,$ where
$$z_1=\ds\frac{1-\sqrt{17}}{4},\quad \ds\frac{1+\sqrt{17}}{4},\quad
z_3=-2.$$
We have
$$z_1=-0.78077\overline{6},\quad
z_2=1.2807\overline{8},\quad z_3=-2,$$
and hence only $z_1$ corresponds to a bound state based
on the criteria given in Theorem~2.5(a), i.e.
among the three $z_j$-values, only $z_1$ lies in the
union of the intervals $(-1,0)$ and $(0,1).$
Next, for $n\ge 1$ we determine $M_n$ defined in (3.14).
The integral in (3.14) can be evaluated in terms of the
residues of $S$ at the two poles $z=0$ and $z=z_1$
inside the unit circle $\bold T.$ From (3.14) and (3.16) we see that
the summation term in (3.14) consists of one term only
and it cancels the contribution from the integral
due to the pole at $z=z_1.$ Thus, the only contribution to
$M_n$ comes from the residues at $z=0$ and we have
$$M_1=-\text{Res}\left[S,0
\right],\quad M_2=-\text{Res}\left[z S,0
\right],\quad M_3=-\text{Res}\left[z^2 S,0
\right].\tag 6.38$$
Using (6.37) in (6.38) we get
$$M_1=-\ds\frac{7}{8},\quad M_2=\ds\frac{1}{4},\quad M_3=\ds\frac{1}{2},
\tag 6.39$$
$$M_n= 0,\qquad n\ge 4.\tag 6.40$$
Using (6.39) and (6.40) let us write the
Marchenko system (3.12)
by retaining only the nonzero terms. We get
$$\cases K_{01}+M_1+K_{01}M_2+K_{02}M_3=0,\\
\stretch
K_{02}+M_2+K_{01}M_3=0,\\
\stretch
K_{03}+M_3=0,\\
\stretch
K_{12}+M_3=0.\endcases\tag 6.41$$
Thus, from the Marchenko system (3.12)
we get $K_{nm}= 0$ for $0\le n<m,$ with the
exception of the four nonzero terms
$K_{01},$ $K_{02},$ $K_{03},$ and $K_{12}.$ By solving the
linear system (6.41) we obtain
$$K_{01}=1,\quad K_{02}=-\ds\frac{3}{4},\quad K_{03}=-\ds\frac{1}{2},\quad K_{12}
=-\ds\frac{1}{2}.$$
Finally, we recover the potential via (2.51) and obtain
$$V_1=K_{01}-K_{12},\quad V_2=K_{12}-K_{23},$$
which yields
$$V_1=\ds\frac{3}{2},\quad V_2=-\ds\frac{1}{2}.$$

In the seventh example below we illustrate
the Gel'fand-Levitan method presented in Section~4.

\noindent {\bf Example 6.7} In this example we illustrate the use of
the Gel'fand-Levitan method to recover the potential from the transmission
eigenvalues. Assume that $b=2$  and $V_2\ne 0.$ Hence, $V_n= 0$ for
$n\ge 3.$ We have $2(b-1)$ equal to $2,$ and hence
there are two transmission eigenvalues
as asserted in Theorem~5.2. Let us assume that they are
given by
$$\lambda_1=-1,\quad \lambda_2=4.
\tag 6.42$$
 From (5.5) and (6.42) we see that
the quantity $E$ defined in (5.7) is given by
$$E=\lambda^2-3\lambda-4.\tag 6.43$$
Using (2.8) in (6.43) we express $E$ in terms of
the spectral parameter $z$ and obtain
$$(z-z^{-1})\,E=z^3-z^2-5z+5z^{-1}+z^{-2}-z^{-3}.\tag 6.44$$
 From (6.44), with the help of (5.10)
we extract $(f_0-1)/V_2$ by taking only the positive powers of $z$ and get
$$\ds\frac{f_0-1}{V_2}=  z^3-z^2-5z.
\tag 6.45$$
As seen from Theorem~5.3(b), the coefficient of $z$ on the right-hand
side in (6.45) gives us $K_{01}/V_2,$ and hence we have
$$\ds\frac{K_{01}}{V_2}=-5.\tag 6.46$$
Using (6.42) and (6.46) in (5.9), we write (5.9) as
$$3=4+\left(-5-1\right)\,V_2,$$
and hence recover $V_2$ as
$$V_2=\ds\frac{1}{6}.\tag 6.47$$
Note that the first equality in (2.48) indicates that
the coefficient of $z^2$ in (6.45) must be $V_1,$ i.e. we have
$$V_1=-1.$$
However, if $b$ were a large integer, we would not be able to obtain
the ordered set
$\{V_1,V_2,\dots,V_{b-1}\}$ of potential values in such a straightforward manner. Instead,
we will use the Gel'fand-Levitan procedure to
illustrate the recovery of $V_1$ because
that is a fundamental procedure to obtain the potential.
Using (6.47) in (6.45) we obtain $f_0$ as
$$f_0=\ds\frac{1}{6}\, z^3-\ds\frac{1}{6}\, z^2+\ds\frac{5}{6}\,z+1.
\tag 6.48$$
 From (2.46) and (2.49) we see that the quantity $g_0$ is obtained
by replacing $z$ with $z^{-1}$ in (6.48) and hence we have
$$g_0=\ds\frac{1}{6}\, z^{-3}-\ds\frac{1}{6}\, z^{-2}+\ds\frac{5}{6}\,z^{-1}+1.
\tag 6.49$$
We determine the zeros of $f_0$ inside the unit circle $|z|=1.$ Note
that $f_0$ has three zeros at $z=z_s$ for $s=1,2,3,$ where
$$z_1=\ds\frac{-1-\sqrt{13}}{2},\quad \ds\frac{-1+\sqrt{13}}{2},\quad
z_3=2.\tag 6.50$$
We have
$$z_1=-2.3027\overline{8},\quad
z_2=1.3027\overline{8},\quad z_3=2,$$
and hence among the three $z_j$-values none lie in the
union of the real intervals $(-1,0)$ and $(0,1).$ Thus,
by Theorem~2.5(a) there are no bound states.
Let us now evaluate $G_{nm}$ given in (4.15). Since there
are no bound states, from (4.3) and (4.12) we obtain
$$d(\rho-\overset{\circ}\to\rho)=\ds\frac{1}{2\pi i}\,(z-z^{-1})^2\,
\left(1-\ds\frac{1}{f_0\,g_0}\right)\,\ds\frac{dz}{z},
\qquad z\in \overline{\bold T^+},\tag 6.51$$
which is extended in a straightforward manner from $\bold T^+$
to the whole circle $\bold T.$
Using (2.21) and (6.51) in (4.15) we get
$$G_{nm}=\ds\frac{1}{2}\,\ds\frac{1}{2\pi i } \oint  \ds\frac{dz}{z}\,
\left(1-\ds\frac{1}{f_0\,g_0}\right)\,\left( z^n-z^{-n}\right)
\,\left( z^m-z^{-m}\right),\tag 6.52$$
where the factor $1/2$ in front of the integral in (6.52) is due to the fact
that we now integrate along the entire unit circle $\bold T$
in the positive direction
instead of integrating along the upper semicircle $\bold T^+.$
We remark the symmetry $G_{nm}=G_{mn},$ as seen from (6.52).
The integral in (6.52) can be evaluated by using the residues
at the poles of the integrand inside the unit circle. Such poles occur
at $z=0$ and $z=1/z_j$ for $j=1,2,3,$ where the $z_j$-values are given in (6.50).
We only need to determine $V_1$ and $V_2$ via the Gel'fand-Levitan method,
and hence
it is sufficient to determine $G_{nm}$ only for
some small values of $n$ and $m.$
Since we will use (2.84) for $n=1$ and $n=2,$ from (2.84) we see that
it is enough for us to set up
the Gel'fand-Levitan system only for $1\le m<n\le 3.$
Using (6.48) and (6.49) in (6.52) we obtain
$$G_{11}=0,\quad G_{21}=1,\quad G_{22}=1, \quad G_{31}=1,\quad G_{32}=\ds\frac{11}{6}.
\tag 6.53$$
We write the Gel'fand-Levitan system (4.17) for $1\le m<n\le 3$ as
$$\cases A_{21}+G_{21}+A_{21}\,G_{11}=0,\\
\stretch
A_{31}+G_{31}+A_{31}G_{11}+A_{32}\,G_{21}=0,\\
\stretch
A_{32}+G_{32}+A_{31}G_{12}+A_{32}\,G_{22}=0.\endcases\tag 6.54$$
Using (6.53) in (6.54) we get
$$\cases A_{21}+1=0,\\
\stretch
A_{31}+1+A_{32}=0,\\
\stretch
A_{32}+\ds\frac{11}{6}+A_{31}+A_{32}=0,\endcases$$
which is uniquely solvable and yielding
$$A_{21}=-1,\quad A_{31}=-\ds\frac{1}{6},\quad A_{32}=-\ds\frac{5}{6}.
\tag 6.55$$
We recall that $A_{10}= 0$ as a result of
(2.81). From (2.84) we have
$$V_1=A_{21}-A_{10},\quad V_2=A_{32}-A_{21}
,\tag 6.56$$
and hence using (6.55) and (6.56) we get
$$V_1=-1,\quad V_2=\ds\frac{1}{6}.$$

\vskip 5 pt

\noindent {\bf Acknowledgments.} This work has been partially supported by the grant DMS-1347475 from the National Science
Foundation.

\vskip 5 pt

\noindent {\bf{References}}

\item{[1]}
T. Aktosun, D. Gintides, and V. G. Papanicolaou, {\it The uniqueness in the inverse problem for transmission eigenvalues for the spherically symmetric variable-speed wave equation,}
Inverse Problems {\bf 27}, 115004 (2011).

\item{[2]}
T. Aktosun, M. Klaus, and R. Weder, {\it Small-energy analysis for the self-adjoint matrix Schr\"odinger operator on the half line,}
J. Math. Phys. {\bf 52}, 102101 (2011).

\item{[3]}
T. Aktosun and V. G. Papanicolaou, {\it Reconstruction of the wave speed from transmission eigenvalues for the spherically symmetric variable-speed wave equation,}
Inverse Problems {\bf 29}, 065007 (2013).

\item{[4]}
T. Aktosun and V. G. Papanicolaou, {\it  Transmission eigenvalues for the self-adjoint Schr\"odinger operator on the half line,}
Inverse Problems {\bf 30}, 075001 (2014).

\item{[5]} T. Aktosun and R. Weder, {\it
Inverse spectral-scattering problem with two sets of discrete spectra for the radial Schr\"odinger equation,} Inverse Problems {\bf 22}, 89--114
(2006).

\item{[6]} C. M. Bender and S. A. Orszag, {\it Advanced mathematical methods for scientists and engineers,} Springer, New York, 1999.

\item{[7]}
F. Cakoni, D. Colton, and H. Haddar,
{\it On the determination of Dirichlet or transmission eigenvalues from far field data,}
C. R. Math. Acad. Sci. Paris {\bf 348}, 379--383 (2010).

\item{[8]}
F. Cakoni, D. Colton, and P. Monk,
{\it
On the use of transmission eigenvalues to estimate the index of refraction from far field data,} Inverse Problems {\bf 23}, 507--522 (2007).

\item{[9]}
K. M. Case and M. Kac,
{\it A discrete version of the inverse scattering problem,}
J. Math. Phys. {\bf 14}, 594--603 (1973).

\item{[10]} K. Chadan and P. C. Sabatier, {\it Inverse problems in
quantum scattering theory,} 2nd ed., Springer, New York, 1989.

\item{[11]} E. A. Coddington and N. Levinson, {\it Theory of ordinary
differential equations,} McGraw-Hill, New York, 1955.

\item{[12]}
D. Colton and R. Kress,
{\it Inverse acoustic and electromagnetic scattering theory,}
2nd ed.,  Springer,
New York, 1998.

\item{[13]}
D. Colton and Y. J. Leung,
{\it
Complex eigenvalues and the inverse
spectral problem for transmission eigenvalues,} Inverse Problems {\bf 29}, 104008 (2013).

\item{[14]}
D. Colton and P. Monk,
{\it The inverse scattering problem for time-harmonic acoustic waves in an inhomogeneous medium,}
Quart. J. Mech. Appl. Math. {\bf 41}, 97--125 (1988).

\item{[15]}
D. Colton, L. P\"aiv\"arinta, and J. Sylvester,
{\it The interior transmission problem,}
Inverse Probl. Imaging {\bf 1}, 13--28 (2007).

\item{[16]} I. M. Gel'fand and B. M. Levitan, {\it On the
determination of a differential equation from its spectral
function,} Amer. Math. Soc. Transl. {\bf 1} (ser. 2),
253--304 (1955).

\item{[17]} B. M. Levitan, {\it
Inverse Sturm Liouville Problems,} VNU Science Press, Utrecht, 1987.

\item{[18]} V. A. Marchenko,
{\it Sturm-Liouville operators and applications,}
Birkh\"auser, Basel, 1986.

\item{[19]}
J. R. McLaughlin and P. L. Polyakov,
{\it On the
uniqueness of a spherically symmetric speed of sound from transmission
eigenvalues,} J. Differential Equations {\bf 107}, 351--382 (1994).

\item{[20]}
J. R. McLaughlin, P. L. Polyakov, and P. E. Sacks,
{\it Reconstruction of a spherically symmetric speed of sound,}
SIAM J. Appl. Math. {\bf 54}, 1203--1223 (1994).

\item{[21]}
J. R. McLaughlin, P. E. Sacks, and M. Somasundaram,
{\it Inverse scattering in acoustic media using interior transmission eigenvalues,}
in: G. Chavent, G. Papanicolaou, P. Sacks, and W. Symes (eds.),
{\it Inverse problems in wave propagation,} Springer,
New York, 1997, pp. 357--374.

\item{[22]}
V. G. Papanicolaou and A. V. Doumas, {\it  On the discrete one-dimensional inverse transmission eigenvalue problem,}
Inverse Problems {\bf 27}, 015004 (2011).

\item{[23]} M. Reed and B. Simon,
{\it Functional analysis,}
Academic Press, New York, 1980.

\end